\makeatletter
\newcommand{\dontusepackage}[2][]{%
  \@namedef{ver@#2.sty}{9999/12/31}%
  \@namedef{opt@#2.sty}{#1}}
\makeatother
\dontusepackage{subfigure}

\documentclass[]{article}

\usepackage{lmodern}
\usepackage{amssymb,amsmath}
\usepackage{ifxetex,ifluatex}
\usepackage[usenames,dvipsnames]{color}
\usepackage{fixltx2e} 
\ifnum 0\ifxetex 1\fi\ifluatex 1\fi=0 
  \usepackage[T1]{fontenc}
  \usepackage[utf8]{inputenc}
\else 
  \ifxetex
    \usepackage{mathspec}
    \usepackage{xltxtra,xunicode}
  \else
    \usepackage{fontspec}
  \fi
  \defaultfontfeatures{Mapping=tex-text,Scale=MatchLowercase}
  
\fi
\IfFileExists{upquote.sty}{\usepackage{upquote}}{}
\IfFileExists{microtype.sty}{%
\usepackage{microtype}
\UseMicrotypeSet[protrusion]{basicmath} 
}{}
\usepackage[margin=1.0in,bottom=1.5in]{geometry}
\usepackage[]{natbib}
\usepackage{listings}
\usepackage{graphicx}
\makeatletter
\def\maxwidth{\ifdim\Gin@nat@width>\linewidth\linewidth\else\Gin@nat@width\fi}
\def\maxheight{\ifdim\Gin@nat@height>\textheight\textheight\else\Gin@nat@height\fi}
\makeatother
\setkeys{Gin}{width=\maxwidth,height=\maxheight,keepaspectratio}
\usepackage{caption}
\usepackage{float}
\usepackage{subfig}
\usepackage{algorithm} 
\let\scholmdAlgorithm\algorithm
\let\endscholmdAlgorithm\endalgorithm
\let\algorithm\relax \let\endalgorithm\relax

{
 \catcode`\*=11\relax
 \global\let\scholmdAlgorithm*\algorithm*
 \global\let\endscholmdAlgorithm*\endalgorithm*
 \global\let\algorithm*\relax 
 \global\let\endalgorithm*\relax
}
\ifxetex
  \usepackage[setpagesize=false, 
              unicode=false, 
              xetex]{hyperref}
\else
  \usepackage[unicode=true]{hyperref}
\fi
\hypersetup{breaklinks=true,
            bookmarks=true,
            pdfauthor={},
            pdftitle={Lift and Relax for PDE-constrained inverse problems in seismic imaging},
            colorlinks=true,
            citecolor=black,
            urlcolor=blue,
            linkcolor=black,
            pdfborder={0 0 0}}
\newcommand{\km}{\mathrm{km}}

\newcommand{\mvec}{\mathbf{m}}
\newcommand{\mvect}{\mathbf{\tilde{m}}}
\newcommand{\uvec}{\mathbf{u}}
\newcommand{\uvect}{\mathbf{\tilde{u}}}
\newcommand{\dvec}{\mathbf{d}}
\newcommand{\rvec}{\mathbf{r}}

\newcommand{\qvec}{\mathbf{q}}
\newcommand{\pvec}{\mathbf{p}}
\newcommand{\svec}{\mathbf{s}}

\newcommand{\nsrc}{n_{\text{s}}}
\newcommand{\nrec}{n_{\text{r}}}
\newcommand{\niter}{n_{\text{it}}}
\newcommand{\ngrid}{n_{\text{g}}}

\newcommand{\nfreq}{n_{\text{f}}} 
\newcommand{\nus}{n_{\text{u}}}

\newcommand{\Pmat}{\mathbf{P}}
\newcommand{\Xmat}{\mathbf{X}}

\newcommand{\Amat}{\mathbf{A}}
\newcommand{\Amatt}{\tilde{\mathbf{A}}}

\newcommand{\Rmat}{\mathbf{R}}

\newcommand{\CC}{\mathbb{C}}
\newcommand{\RR}{\mathbb{R}}
\newcommand{\Smat}{\mathbf{S}}
\newcommand{\Tmat}{\mathbf{T}}
\newcommand{\Smatt}{\mathbf{\tilde{S}}}
\newcommand{\mHz}{\mathrm{Hz}}

\newcommand{\fpen}{f_\text{p}(\mvec,\uvec)}
\newcommand{\fpent}{f_{\text{p}_{2}}(\mvect,\uvect,\alpha)}
\newcommand{\fpentt}{f_{\text{p}_{2}}(\mvect,\uvect,\theta)}
\newcommand{\fpentts}{f_{\text{p}_{2}}(\mvect^{\star},\uvect,\theta^{\star})}
\newcommand{\fpenttu}{f_{\text{p}_{2}}(\mvect,\uvect^{\star}(\mvect, \theta),\theta)}
\newcommand{\fpentb}{\overline{f}_{\text{p}_{2}}(\mvect, \theta)}
\newcommand{\fpentbk}{\overline{f}^{(k)}_{\text{p}_{2}}(\mvect^{(k)}, \theta^{(k)})}
\newcommand{\fpenx}{f_\text{px}(\Xmat)}
\newcommand{\fpenr}{f_\text{pr}(\Rmat)}
\newcommand{\fred}{f_\text{r}(\mvec)}
\newcommand{\diag}{\text{diag}}
\newcommand{\full}{f_\text{f}(\mvec,\uvec)}
\newcommand{\mkm}{\mathrm{km}}

\title{Lift and Relax for PDE-constrained inverse problems in seismic imaging}
\author{Zhilong Fang\textsuperscript{1,2}, Laurent
Demanet\textsuperscript{1,2}\\\textsuperscript{1}Department of
Mathematics, Massachusetts Institute of
Technology\\\textsuperscript{2}Earth Resource Laboratory, Massachusetts
Institute of Technology}
\date{}

\begin{document}
\maketitle

\section{Abstract}\label{abstract}

We present Lift and Relax for Waveform Inversion (LRWI), an approach
that mitigates the local minima issue in seismic full waveform inversion
(FWI) via a combination of two convexification techniques. The first
technique (Lift) extends the set of variables in the optimization
problem to products of those variables, arranged as a moment matrix.
This algebraic idea is a celebrated way to replace a hard polynomial
optimization problem by a semidefinite programming approximation.
Concretely, both the model and the wavefield are lifted from vectors to
rank-2 matrices. The second technique (Relax) invites to consider the
wave equation, not as a hard constraint, but as a soft constraint to be
satisfied only approximately -- a technique known as wavefield
reconstruction inversion (WRI). WRI weakens wave-equation constraints by
introducing wave-equation misfits as a weighted penalty term in the
objective function. The relaxed penalty formulation enables balancing
the data and wave-equation misfits by tuning a penalty parameter.
Together, ``Lift'' and ``Relax'' help reformulate the inverse problem as
a set of constraints on a rank-2 moment matrix in a higher dimensional
space. Such a lifting strategy permits a good data and wave-equation fit
throughout the inversion process, while leaving the numerical rank of
the rank-2 moment matrix to be minimized down to one. Numerical examples
indicate that compared to FWI and WRI, LRWI can conduct successful
inversions using an initial model that would be considered too poor, and
data with a starting frequency that would be considered too high, for
either method in isolation. Specifically, LRWI increases the acceptable
starting frequency from 1.0 Hz and 0.5 Hz to 2.0 Hz and 2.5 for the
Marmousi model and the Overthrust model, respectively, in the cases of a
linear gradient starting model.

\section{Key words.}\label{key-words.}

Inverse problem, full waveform inversion, lift, relax, PDE

\section{AMS subject classification.}\label{ams-subject-classification.}

86A22, 35R30

\section{Introduction}\label{introduction}

Seismic imaging is the primary means for Earth scientists and
geophysicists to explore and study Earth's deep interior, where direct
observations are infeasible. Its applications range from studies of
Earth's core, thousands of kilometers below the surface, to detailed
images of shallow crustal structures for locating petroleum deposits.
During the last thirty years, with the advancements in high-performance
computing and the development of wide-aperture and dense data
acquisition, seismic imaging techniques have been upgraded from simple
and low-resolution ray-based methods to complicated and high-resolution
wave-equation-based methods. Especially, during the last two decades,
full waveform inversion (FWI)
\citep{Tarantola1982FWI, pratt1999seismic, VirieuxOverview2009} has
become one of the most important approaches because of its potential
capability in creating high-resolution subsurface images through the
usage of all kinds of waves in the data.

Conventional FWI seeks a subsurface velocity model that can minimize the
difference between its predicted data and the observed data in a
least-squares sense. A well-known problem associated with conventional
FWI is that it suffers from local minima in the objective function
caused by the so-called ``cycle-skipping'' issues. More specifically, if
the initial model does not generate predicted data within half a
wavelength of the observed data, iterative optimization approaches may
stagnate at physically meaningless solutions with a high probability. In
order to conduct a successful inversion, conventional FWI needs a good
initial model that is kinematically accurate at the longest data
wavelengths and data containing enough low frequencies and long offsets
\citep{VirieuxOverview2009, vigh20093d, warner2013full}. Research aimed
at mitigating the ``cycle-skipping'' issue mainly focuses on different
misfit functions
\citep{cara1987waveform, van2010correlation, wu2013ultra, engquist2014application, warner2016adaptive, yang2018application},
expanding the search space
\citep{vanleeuwen2015IPpmp, huang2017full, fang2018source, fang2018uncertainty},
and the integration with the advanced approach of migration velocity
analysis \citep{symes2008migration, li2014wave}.

We propose a two-pronged Lift and Relax waveform inversion (LRWI)
approach to mitigating the local minima problem in this paper. The
proposed approach consists of two relaxation strategies that expand the
search space. The ``Relax'' strategy is based on the so-called approach
wavefield reconstruction inversion (WRI)
\citep{vanleeuwen2015IPpmp, fang2018source}. WRI first introduces
wavefields as additional unknown variables, and then weakens the partial
differential equation (PDE) constraints used in conventional FWI by
treating the PDE misfit as a weighted penalty term in the objective
function. Through tuning the penalty parameter, the resulting approach
does not enforce the PDE constraints at each iteration and arguably
yields a less non-linear problem in the model parameter. The ``Lift''
strategy follows the early work in \citet{cosse2015short}~that borrows
ideas from recent developments in the semidefinite relaxation for
polynomial equations to mitigate non-convexity
\citep{lasserre2001global, laurent2009sums}. We lift both unknown
wavefields and model parameters from 1D vectors to rank-2 matrices, and
reformulate the WRI problem as a set of constraints on a rank-2 moment
matrix in a higher dimensional space. Such a lifting strategy permits a
good data and wave-equation fit throughout the inversion process, while
leaving the numerical rank of the moment matrix to be the quantity to
minimize -- so that this matrix aims to be a rank one matrix at
convergence eventually.

Compared to conventional FWI, the proposed LRWI approach has three major
advantages. First, the computation of the gradients does not require
adjoint or reverse-time wavefields. Secondly, the ``Relax'' and ``Lift''
strategies enable us to fit both data misfit and PDE misfit even with
poor models. Thirdly, the rank-2 formulation provides us with the
potential to utilize information from the two components in the rank-2
model matrix simultaneously. The last two properties, in conjunction
with the expanded search space, may result in an optimization
formulation that is less prone to local minima. We present numerical
examples on both Marmousi and Overthrust models to illustrate the
feasibility and advantages of the proposed approach.

The paper is organized as follows. First, we review the basic conception
and formulation of conventional FWI. Next, we derive the formulation for
the proposed rank-2 LRWI. Then, we derive all the necessary components
for the efficient optimization strategy in detail. Finally, we present
numerical examples on Marmousi and Overthrust models to illustrate the
feasibility and advantages of LRWI and conclude the paper with a
detailed discussion.

\section{Methodology}\label{methodology}

Given a seismic data set
$\dvec \in \mathrm{R}^{\nsrc\times\nrec\times\nfreq}$ with $\nsrc$
sources, $\nrec$ receivers, and $\nfreq$ frequencies, FWI aims to
reconstruct the discretized $\ngrid$-dimensional squared slowness model
$\mvec$ from $\dvec$ by solving the following PDE-constrained
optimization problem:
\begin{equation}
\begin{aligned}
    &\min_{\mvec, \uvec}\full=\frac{1}{2}\sum_{i,j}^{\nsrc,\nfreq}\|\Pmat\uvec_{i,j}-\dvec_{i,j}\|^2_{2},\\
    &\text{subject to} \quad (\Delta+\omega_{j}^{2}\mvec)\uvec_{i,j} = \qvec_{i,j},
\end{aligned}
\label{FWI}
\end{equation}
 where the operator $\Pmat$ projects the wavefield $\uvec_{i,j}$
corresponding to the $i^{\text{th}}$ source $\qvec_{i,j}$ with frequency
$\omega_{j}$ onto the receiver locations. The operator $\Delta$
represents the Laplacian operator, and the equation
$(\Delta+\omega_{j}^{2}\mvec)\uvec_{i,j} = \qvec_{i,j}$ is known as the
Helmholtz equation.

The optimization problem in Equation~\ref{FWI} requires a solution in
$\RR^{\ngrid} \times \CC^{\nus}$ with
$\nus = \nsrc\times\nfreq\times\ngrid$, which is infeasible for most
practical applications because we cannot afford to store all the unknown
variables. To reduce the dimensionality of the search space, the
conventional adjoint-state method \citep{VirieuxOverview2009} eliminates
the PDE constraint
$(\Delta+\omega_{j}^{2}\mvec)\uvec_{i,j} = \qvec_{i,j}$ through solving
the PDE straightforwardly, yielding the following reduced problem:
\begin{equation}
\begin{aligned}
    \min_{\mvec}\fred&=\frac{1}{2}\sum_{i,j}^{\nsrc,\nfreq}\|\Pmat\Amat_j(\mvec)^{-1}\qvec_{i,j}-\dvec_{i,j}\|^2_{2}, \\ &\text{with}\,\, \Amat_j(\mvec) = \Delta+\omega_{j}^{2}\mvec,
\end{aligned}
\label{FWI2}
\end{equation}
 whose search space is $\RR^{\ngrid}$. Although the dimensionality of
the search space reduces from $\nus+\ngrid$ to $\ngrid$, the trade-off
lies in the fact that the inversion of the Helmholtz matrix introduces a
very strong nonlinearity into the problem, yielding an objective
function $\fred$ with many local minima.

\subsection{\texorpdfstring{WRI with a rank-\emph{r}
relaxation}{WRI with a rank-r relaxation}}\label{wri-with-a-rank-r-relaxation}

In this work, we aim to mitigate the local minima issue of conventional
FWI by proposing a Lift and Relax formulation in the rank-$r$ case. To
simplify the notation, we will omit the dependence of the variables on
the source and frequency indexes $i$ and $j$ from now on.

We first follow \citet{vanleeuwen2015IPpmp} and relax the PDE constraint
in Equation~\ref{FWI} by considering the PDE misfit as a weighted
penalty term as follows:
\begin{equation}
    \min_{\mvec, \uvec}\fpen=\frac{1}{2}\|\Pmat\uvec-\dvec\|^2_{2}+\frac{\lambda}{2}\|(\Delta+\omega^{2}\mvec)\uvec - \qvec\|^2_{2}.
\label{WRI}
\end{equation}
 The penalty parameter $\lambda$ enables us to balance the PDE and data
misfits and provides the freedom to design a search path in the enlarged
space that can potentially bypass the local minima in the objective
function of conventional FWI.

Following the PDE relaxation, we introduce an additional rank-$r$
relaxation to expand the search space into a higher dimension space,
which is motivated from the following matrix expression of the unknown
parameters $\mvec$ and $\uvec$:
\begin{equation}
\begin{aligned}
    \Xmat = \begin{bmatrix}
    \Xmat_{11} & \Xmat_{12} & \Xmat_{13} \\
    \Xmat_{21} & \Xmat_{22} & \Xmat_{23} \\
    \Xmat_{31} & \Xmat_{32} & \Xmat_{33} 
    \end{bmatrix}=[1, \mvec^{\top}, \uvec^{\top}]^{\top}[1, \mvec^{\top}, \uvec^{\top}].
\end{aligned}
\label{Lift1}
\end{equation}
 Clearly, the matrix $\Xmat$ is a rank-1 positive semindefinite matrix.
Based on Equation~\ref{Lift1}, we can lift the original WRI problem from
optimizing over vectors $\mvec$ and $\uvec$ to optimizing over the
matrix $\Xmat$. In the course of doing so, the direct correspondence to
$\mvec$ and $\uvec$ in Equation~\ref{Lift1} is not directly imposed, but
the objective in Equation~\ref{WRI} is rewritten with the blocks of
$\Xmat$ serving as proxies for $\mvec$, $\uvec$, and the product
$\mvec\uvec^{\top}$. This yields the following equivalent optimization
problem:
\begin{equation}
\begin{aligned}
    \min_{\Xmat}\fpenx=\frac{1}{2}\|\Pmat\Xmat_{31}-\dvec\|^2_{2}&+\frac{\lambda}{2}\|\Delta\Xmat_{31}+\omega^{2}\diag(\Xmat_{32}) - \qvec\|^2_{2},\\
    \text{subject to}\quad \Xmat_{11} &= 1, \\
    \Xmat &\succeq 0, \\
    \text{rank}(\Xmat) &= 1.
\end{aligned}
\label{WRIX}
\end{equation}
 The ``Lift'' relaxation then consists in dropping the rank-1
constraint.

The new objective function $\fpenx$ is a quadratic function with respect
to the matrix $\Xmat$, which is much simpler than the original FWI and
WRI objective functions. Since $\Xmat \in \CC^{(\nus + \ngrid + 1)^2}$,
we are not able to optimize over $\Xmat$ directly for large-scale
realistic applications. Nonetheless, as stated by
\citet{cosse2015short}, it is possible for us to obtain a
computationally feasible formulation with a reasonable storage
requirement by introducing a rank-$r$ factorization $\Rmat\Rmat^{\top}$
for the matrix $\Xmat$:
\begin{equation}
\begin{aligned}
    \min_{\Rmat}\fpenr=\frac{1}{2}\|\Pmat(\Rmat_{3}\Rmat_{1}^{\top})-\dvec\|^2_{2}&+\frac{\lambda}{2}\|\Delta(\Rmat_{3}\Rmat_{1}^{\top})+\omega^{2}\diag(\Rmat_{3}\Rmat_{2}^{\top}) - \qvec\|^2_{2},\\
    \text{subject to}\quad \Rmat_{1}\Rmat_{1}^{\top} &= 1, 
\end{aligned}
\label{WRIR}
\end{equation}
 where
$\Rmat = (\Rmat_{1}^{\top}, \Rmat_{2}^{\top}, \Rmat_{3}^{\top})^{\top}$
with
$\Rmat_{1} = [\alpha_1,..., \alpha_r]\in \RR^{1\times r}, \, \Rmat_{2}=[\mvect_{1},...,\mvect_{r}]\in \RR^{\ngrid\times r}, \text{and}\, \Rmat_{3}=[\uvect^1,..,\uvect^{r}]\in\CC^{\nus\times r}$.
This block representation of $\Rmat$ leads to a representation of
$\Xmat$ as a sum of rank-1 matrices,
\begin{equation}
\begin{aligned}
    \Xmat\approx \Rmat\Rmat^{\top} = \sum_{l=1}^{r}\begin{bmatrix} 
    \alpha_l^2,       & \alpha_l \mvect_l^{\top} & \alpha_l \uvect_l^{\top} \\
    \alpha_l\mvect_l  & \mvect_l \mvect^{\top}_l  & \mvect_l \uvect_l^{\top} \\
    \alpha_l \uvect_l & \uvect_l \mvect^{\top}_l  & \uvect_l \uvect_l^{\top}
    \end{bmatrix}.
\end{aligned}
\label{RankrRelax}
\end{equation}
 When $r = 1$, the optimization problem in Equation~\ref{WRIR} will
reduce to the original WRI problem in Equation~\ref{WRI}. A larger $r$
yields a stronger relaxation but introduces more computational cost and
storage requirements.

\subsubsection{rank-2 relaxation}\label{rank-2-relaxation}

In this work, we present a rank-2 formulation for the optimization
problem in Equation~\ref{WRIR} to balance the relaxation and
computational costs. When selecting $r = 2$, we have
\begin{equation}
\begin{aligned}
    \mvec           &= \alpha_1 \mvect_1 + \alpha_2 \mvect_2, \\
    \uvec           &= \alpha_1 \uvect_1 + \alpha_2 \uvect_2, \\
    \mvec\odot\uvec &= \mvect_1\odot\uvect_1 + \mvect_2\odot\uvect_2, \\
    1               &= \alpha_1^2 + \alpha_2^2,
\end{aligned}
\label{Relax1}
\end{equation}
 where the operator $\odot$ represents the pointwise multiplication or
the Hadamard product. The rank-2 expression in Equation~\ref{Relax1}
yields the following optimization problem:
\begin{equation}
    \min_{\mvect, \uvect, \alpha}\fpent=\frac{1}{2}\|\sum_{l=1}^{2}\Pmat\alpha_l\uvect_l-\dvec\|^2_{2}+\frac{\lambda}{2}\|\sum_{l=1}^{2}\alpha_l\Delta\uvect_l+\omega^{2}\sum_{l=1}^{2}\mvect_{l}\odot\uvect_{l} - \qvec\|^2_{2}.
\label{WRIrank2}
\end{equation}
 It is easy to verify that this optimization problem has infinite
solutions. Indeed, for any fixed pair of ($\mvect^{\star}$,
$\alpha^{\star}$), the optimal $\uvect^{\star}$ for the objective
function $f_{\text{p}_2}(\mvect^{\star}, \uvect, \alpha^{\star})$ should
satisfy the following equation:
\begin{equation}
\begin{aligned}
    \Smat\uvect^{\star}&=
    \begin{bmatrix} 
    \alpha_1\Pmat,                                  & \alpha_2 \Pmat\\
    \lambda^{\frac{1}{2}}(\alpha_1\Delta+\omega^2\mvect_1)  & \lambda^{\frac{1}{2}}(\alpha_2\Delta+\omega^2\mvect_2) 
    \end{bmatrix}\ 
    \begin{bmatrix}
    \uvect^{\star}_{1} \\
    \uvect^{\star}_{2}
    \end{bmatrix} = 
    \begin{bmatrix}
    \dvec \\
    \lambda^{\frac{1}{2}}\qvec
    \end{bmatrix}.
\end{aligned}
\label{optWRIrank2}
\end{equation}
 Since the matrix $\Smat$ is an underdetermined
$(\ngrid+\nrec) \times 2\ngrid$ matrix with $\ngrid > \nrec$, the linear
Equation~\ref{optWRIrank2} has infinite solutions for $\uvect^{\star}$.
As a result, there are infinite global minima ($\mvect^{\star}$,
$\uvect^{\star}$, $\alpha^{\star}$)s satisfying
$f_{\text{p}_2}(\mvect^{\star}, \uvect^{\star}, \alpha^{\star}) = 0$.

To mitigate the nonuniqueness issue of optimizing
Equation~\ref{WRIrank2}, we need additional information to regularize
the problem. We notice that the original lifted problem in
Equation~\ref{WRIX} has the constraint of $\text{rank}(\Xmat) = 1$,
which is not involved in the rank-2 formulation. Therefore, to derive
our regularization, we reimpose this information. We do not
straightforwardly require $\text{rank}(\Rmat) = 1$, otherwise it will
downgrade the problem to the rank-1 case, which is the original WRI
problem. Instead, we use another necessary condition for a rank-1 matrix
to introduce a weaker regularization. If the matrix $\Rmat$ is a rank-1
matrix, then its three components $\Rmat_{1} = [\alpha_1, \alpha_2]$,
$\Rmat_{2} = [\mvect_1, \mvect_2]$ and
$\Rmat_{3} = [\uvect_1, \uvect_2]$ should satisfy the following
requirements:

\begin{equation}
\begin{aligned}
    \alpha_1\mvect_2 &= \alpha_2\mvect_1, \\
    \alpha_1\uvect_2 &= \alpha_2\uvect_1, \\
    \mvect_1\odot\uvect_2 &= \mvect_2\odot\uvect_1.
\end{aligned}
\label{optWRIrank3}
\end{equation}

We can use these properties to regularize the problem. In this work,
since we are more interested in $\mvect$ and $\uvect$ than $\alpha$, we
use the third property to introduce an additional regularization to the
optimization problem in Equation~\ref{WRIrank2} as follows:
\begin{equation}
\begin{aligned}
    \min_{\mvect, \uvect, \alpha}\fpent=&\frac{1}{2}\|\sum_{l=1}^{2}\Pmat\alpha_l\uvect_l-\dvec\|^2_{2}\\&+\frac{\lambda}{2}\|\sum_{l=1}^{2}\alpha_l\Delta\uvect_l+\omega^{2}\sum_{l=1}^{2}\mvect_{l}\odot\uvect_{l} - \qvec\|^2_{2}\\&+\frac{\gamma}{2}\|\mvect_1\odot\uvect_2-\mvect_2\odot\uvect_1\|^{2}_{2}, \\
    &\text{subject to} \quad \alpha_1^2 + \alpha_2^2 = 1,
\end{aligned}
\label{WRIrank3}
\end{equation}
 where the penalty parameter $\gamma$ controls the strength of the
rank-1 regularization.

Finally, we can simplify the constrained optimization problem in
Equation~\ref{WRIrank2} to an unconstrained problem by eliminating the
constraint $\alpha_1^2 + \alpha_2^2 = 1$ with a simple polar coordinates
transform:
\begin{equation}
    \alpha_1 = \sin\theta \quad \text{and} \quad \alpha_2 = \cos\theta,
\label{SinCos}
\end{equation}
 yielding the following unconstrained optimization problem:
\begin{equation}
\begin{aligned}
    \min_{\mvect, \uvect, \theta}\fpentt=&\frac{1}{2}\|\Pmat(\sin\theta\uvect_1+\cos\theta\uvect_2)-\dvec\|^2_{2} \\
                                       +&\frac{\lambda}{2}\|\Delta(\sin\theta\uvect_1+\cos\theta\uvect_2)+\omega^{2}\sum_{l=1}^{2}\mvect_{l}\odot\uvect_{l} - \qvec\|^2_{2}\\
                                       +&\frac{\gamma}{2}\|\mvect_1\odot\uvect_2-\mvect_2\odot\uvect_1\|^{2}_{2}.
\end{aligned}
\label{WRIrank4}
\end{equation}

\subsection{Variable projection and optimization
scheme}\label{variable-projection-and-optimization-scheme}

The optimization problem in Equation~\ref{WRIrank4} still faces the
challenge of a large storage requirement. In order to reduce the storage
requirement, we use the variable projection method
\citep{golub2003separable} to project out the wavefields $\uvect$, which
is the main source of the storage cost. For any pair of
($\mvect^{\star}$, $\theta^{\star}$), the objective function $\fpentts$
is quadratic with respect to $\uvect$, whose minimizer has an analytical
solution:

\begin{equation}
\begin{split}
    \uvect^{\star}= (\Smatt^{\top}\Smatt)^{-1}\Smatt^{\top}\begin{bmatrix}
    \dvec \\
    \lambda^{\frac{1}{2}}\qvec \\
    0
    \end{bmatrix},
\end{split}
\label{AnaSolution}
\end{equation}

with
\begin{equation}
\begin{aligned}
    &\Smatt = \begin{bmatrix} 
    \sin\theta\Pmat                            & \cos\theta \Pmat\\
    \lambda^{\frac{1}{2}}\Amatt(\mvect_1)      & \lambda^{\frac{1}{2}}\Amatt(\mvect_2) \\
    \gamma^{\frac{1}{2}}\text{diag}(\mvect_2)  & -\gamma^{\frac{1}{2}}\text{diag}(\mvect_1)
    \end{bmatrix}, \\
    \quad\text{with} \quad \Amatt(\mvect_{1}) = &\sin\theta\Delta+\omega^2\mvect_1, \quad \text{and} \quad \Amatt(\mvect_{2})= \cos\theta\Delta+\omega^2\mvect_2.    
\end{aligned}
\label{Sexp}
\end{equation}

Replacing the variable $\uvect$ in Equation~\ref{WRIrank4} by the
optimal solution $\uvect^{\star}(\mvect, \theta)$, we obtain a reduced
objective function $\fpentb = \fpenttu$. We can use the chain rule to
compute the derivatives of $\nabla_{\mvect}{\fpentb}$ and
$\nabla_{\theta}{\fpentb}$ as follows:
\begin{equation}
\begin{aligned}
    \nabla_{\mvect}{\fpentb} &= \nabla_{\mvect}{\fpenttu}\\ &=\nabla_{\mvect}{\fpent}\rvert_{\uvect=\uvect^{\star}} + \nabla_{\uvect}{\fpent}\rvert_{\uvect=\uvect^{\star}}\nabla_{\mvect}\uvect, \\
    \nabla_{\theta}{\fpentb} &= \nabla_{\theta}{\fpenttu} \\&=\nabla_{\theta}{\fpent}\rvert_{\uvect=\uvect^{\star}} + \nabla_{\uvect}{\fpent}\rvert_{\uvect=\uvect^{\star}}\nabla_{\theta}\uvect.
\end{aligned}
\label{Derivative}
\end{equation}
 The most important property of the variable projection method lies in
the fact that $\uvect^{\star}$ minimizes the objective function $\fpent$
for fixed $(\mvect, \theta)$, satisfying the condition
$\nabla_{\uvect}{\fpent}\rvert_{\uvect=\uvect^{\star}}=0$. Therefore, we
can drop out the complicated terms
$\nabla_{\uvect}{\fpent}\rvert_{\uvect=\uvect^{\star}}\nabla_{\mvect}\uvect$
and
$\nabla_{\uvect}{\fpent}\rvert_{\uvect=\uvect^{\star}}\nabla_{\theta}\uvect$
in the expressions of $\nabla_{\mvect}{\fpentb}$ and
$\nabla_{\theta}{\fpentb}$, and simplify them as follows:
\begin{equation}
\begin{split}
    \nabla_{\mvect}{\fpentb} = \nabla_{\mvect}{\fpent}\rvert_{\uvect=\uvect^{\star}}, \\
    \nabla_{\theta}{\fpentb} = \nabla_{\theta}{\fpent}\rvert_{\uvect=\uvect^{\star}}.
\end{split}
\label{Derivative2}
\end{equation}
 Following Equation~\ref{Derivative2}, the expressions for
$\nabla_{\mvect}{\fpentb}$ and $\nabla_{\theta}{\fpentb}$ can be derived
as follows:
\begin{equation}
\begin{aligned}
    \nabla_{\mvect}{\fpentb} = &\begin{bmatrix} 
                                \nabla_{\mvect_1}{\fpentb} \\ 
                                \nabla_{\mvect_2}{\fpentb} \end{bmatrix} \\
                             = &\begin{bmatrix} 
                                \lambda(\omega^2\text{diag}(\uvect^{\star}_1))^{\top}\pvec+\gamma(\text{diag}(\uvect^{\star}_2))^{\top}\svec \\ 
                                \lambda(\omega^2\text{diag}(\uvect^{\star}_2))^{\top}\pvec-\gamma(\text{diag}(\uvect^{\star}_1))^{\top}\svec \end{bmatrix},\\
    \nabla_{\theta}{\fpentb} = & \cos\theta[(\Pmat\uvect^{\star}_1)^{\top}\rvec + \lambda(\Delta\uvect^{\star}_{1})^{\top}\pvec] \\
                             - & \sin\theta[(\Pmat\uvect^{\star}_2)^{\top}\rvec + \lambda(\Delta\uvect^{\star}_{2})^{\top}\pvec],
\end{aligned}
\label{Gradients}
\end{equation}
 where
\begin{equation}
\begin{aligned}
    \pvec & = \Delta(\sin\theta\uvect^{\star}_1+\cos\theta\uvect^{\star}_2)+\omega^{2}\sum_{l=1}^{2}\mvect_{l}\odot\uvect^{\star}_{l} - \qvec, \\
    \svec & = \mvect_1\odot\uvect^{\star}_2-\mvect_2\odot\uvect^{\star}_1, \\
    \rvec & = \Pmat(\sin\theta\uvect^{\star}_1+\cos\theta\uvect^{\star}_2)-\dvec.
\end{aligned}
\label{pands}
\end{equation}
 Once obtained $\uvect^{\star}$, Equations~\ref{Gradients}
and~\ref{pands} imply that the computation of the gradients
$\nabla_{\mvect}{\fpentb}$ and $\nabla_{\theta}{\fpentb}$ only involves
simple and cheap matrix-vector multiplications and does not involve any
additional computationally intensive matrix inverses. Compared to the
conventional adjoint-state method that requires to invert an additional
adjoint Helmholtz matrix to obtain the gradient, the proposed method
reduces computational cost for computing the gradient.

With the derivatives in Equation~\ref{Gradients}, we can use
optimization algorithms like gradient descent and limited-memory
Broyden-Fletcher-Goldfarb-Shanno (l-BFGS) method
\citep{nocedal2006numerical} that only needs the gradient information to
solve the optimization problem. During the optimization, since $\mvect$
and $\theta$ are very different in scale and have very different
sensitivities to the objective function, we propose to update them
alternately. During each iteration, we first conduct an l-BFGS update on
$\mvect$, then we use a gradient descent step to update $\theta$.
Algorithm~\ref{alg:1} illustrates the pseudo code of the two-stage
l-BFGS method.

\begin{scholmdAlgorithm}
~~1.~Initialization~with~$\mvect_1^{(0)}$,~$\mvect_2^{(0)}$~and~$\theta^{(0)}$~\\\hspace*{0.333em}\hspace*{0.333em}2.~\textbf{for}~$k = 1 \rightarrow \niter$\\\hspace*{0.333em}\hspace*{0.333em}3.~~~~~~Compute~$\uvect^{\star(k)}$~by~Equation~\ref{AnaSolution}~\\\hspace*{0.333em}\hspace*{0.333em}4.~~~~~~Compute~$\fpentbk$~and~~$\nabla_{\mvect}{\fpentbk}$~~by~Equations~\ref{Gradients}~~~~\\\hspace*{0.333em}\hspace*{0.333em}5.~~~~~~l-BFGS~step~in~$\mvect^{(k)}$~to~get~$\mvect^{(k+1)}$\\\hspace*{0.333em}\hspace*{0.333em}6.~~~~~~Compute~$\nabla_{\theta}\overline{f}^{(k)}_{\text{p}_{2}}(\mvect^{(k+1)}, \theta^{(k)})$~by~Equations~\ref{Gradients}~~~~\\\hspace*{0.333em}\hspace*{0.333em}7.~~~~~~Gradient~descent~step~in~$\theta^{(k)}$~to~get~$\theta^{(k+1)}$\\\hspace*{0.333em}\hspace*{0.333em}8.~\textbf{end}~~~~~\\\hspace*{0.333em}\hspace*{0.333em}9.~Obtain~$\theta^{\star}$~and~$\mvect^{\star}=(\mvect_{1}^{\star}, \mvect_{2}^{\star})$\\\hspace*{0.333em}10.~Output~$\mvec^{\star} = \sin\theta^{\star}\mvect_{1}^{\star} + \cos\theta^{\star}\mvect_{2}^{\star}$~~~~~~
\caption{Rank-2 LRWI}\label{alg:1}
\end{scholmdAlgorithm}

\subsection{\texorpdfstring{Selection of $\lambda$ and
$\gamma$}{Selection of \textbackslash{}lambda and \textbackslash{}gamma}}\label{selection-of-lambda-and-gamma}

The selection of $\lambda$ and $\gamma$ plays an important role in the
proposed LRWI, because $\lambda$ and $\gamma$ affect the condition
number of the matrix $\Smatt$ in Equation~\ref{AnaSolution} and the
search path. An appropriate selection can produce a search path that
bypasses the local minima of conventional FWI and also speeds up the
optimization procedure. In this work, we propose a two-stage unit-free
strategy to select $\lambda$ and $\gamma$.

We first determine the selection of $\lambda$.
\citet{vanleeuwen2015IPpmp} and \citet{fang2018uncertainty} studied the
selection of $\lambda$ for WRI and proposed a natural scaling for
$\lambda$, i.e.
$\lambda > \mu_{1}(\Amat^{-\top}\Pmat^{\top}\Pmat\Amat^{-1})$ can be
considered large, while
$\lambda < \mu_{1}(\Amat^{-\top}\Pmat^{\top}\Pmat\Amat^{-1})$ can be
considered small, where the matrix $\Amat$ denotes the Helmholtz matrix
parameterized by the current model $\mvec$ and
$\mu_{1}(\Amat^{-\top}\Pmat^{\top}\Pmat\Amat^{-1})$ denotes the largest
eigenvalue of the matrix $\Amat^{-\top}\Pmat^{\top}\Pmat\Amat^{-1}$.
Specifically, when
$\lambda < 10^{-2}\mu_{1}(\Amat^{-\top}\Pmat^{\top}\Pmat\Amat^{-1})$,
the simulated wavefields tend to fit the observed data while leaving a
big misfit for the PDE; when
$\lambda > 10^2\mu_{1}(\Amat^{-\top}\Pmat^{\top}\Pmat\Amat^{-1})$, the
opposite holds. In practice, considering the large computational cost of
calculating $\mu_1$, \citet{vanleeuwen2015IPpmp} suggest using $\mu_1$
parameterized with the initial model $\mvec^{(0)}$ to select the penalty
parameter $\lambda$. Following \citet{vanleeuwen2015IPpmp} and
\citet{fang2018uncertainty}, we select $\lambda$ according to the value
$\mu_{1}(\Amat(\mvec^{(0)})^{-\top}\Pmat^{\top}\Pmat\Amat(\mvec^{(0)})^{-1})$,
where
$\mvec^{(0)} = \sin\theta^{(0)}\mvect^{(0)}+\cos\theta^{(0)}\mvect^{(0)}$.

With $\lambda$ in hand, the selection of $\gamma$ will determine the
condition number of the matrix $\Smatt$. Since both blocks
$\begin{bmatrix}  \sin\theta\Pmat & \cos\theta \Pmat\\  \lambda^{\frac{1}{2}}\Amatt(\mvect_1) & \lambda^{\frac{1}{2}}\Amatt(\mvect_2) \end{bmatrix}$
and
$\begin{bmatrix} \gamma^{\frac{1}{2}}\text{diag}(\mvect_2) & -\gamma^{\frac{1}{2}}\text{diag}(\mvect_1)  \end{bmatrix}$
are underdetermined, either a very large $\gamma$ or a very small
$\gamma$ will lead to a bad conditioned matrix $\Smatt$. Indeed the
matrix $\Smatt^{\top}\Smatt$ in Equation~\ref{AnaSolution} has the
following expression:
\begin{equation}
\begin{aligned}
    \Smatt^{\top}\Smatt &= \begin{bmatrix} 
     \Tmat_{1,1} & \Tmat_{1,2} \\ 
     \Tmat_{2,1} & \Tmat_{2,2} 
    \end{bmatrix}, \quad \text{with}\\  
    \Tmat_{1,1} &= \alpha_1^{2}\Pmat^{\top}\Pmat + \lambda\Amatt(\mvect_1)^{\top}\Amatt(\mvect_1) + \gamma\text{diag}(\mvect_{2}\odot\mvect_{2}), \\
    \Tmat_{1,2} &= \alpha_1\alpha_2\Pmat^{\top}\Pmat + \lambda\Amatt(\mvect_1)^{\top}\Amatt(\mvect_2) - \gamma\text{diag}(\mvect_{1}\odot\mvect_{2}), \\
    \Tmat_{2,1} &= \alpha_1\alpha_2\Pmat^{\top}\Pmat + \lambda\Amatt(\mvect_2)^{\top}\Amatt(\mvect_1) - \gamma\text{diag}(\mvect_{1}\odot\mvect_{2}), \\
    \Tmat_{2,2} &= \alpha_2^{2}\Pmat^{\top}\Pmat + \lambda\Amatt(\mvect_2)^{\top}\Amatt(\mvect_2) + \gamma\text{diag}(\mvect_{1}\odot\mvect_{1}).
\end{aligned}
\label{StS}
\end{equation}
 Equation~\ref{StS} motivates us to derive the scaling of $\gamma$ by
comparing $\gamma\mvect_{i}\odot\mvect_{j}$ with the diagonal part of
the matrices
$\Tmat(\lambda) = \{ \Tmat_{i,j}=\lambda\Amatt(\mvect_{i})^{\top}\Amatt(\mvect_{j})+\alpha_{i}\alpha_{j}\Pmat^{\top}\Pmat\}_{1\leq i,j \leq 2}$.
A natural scaling for $\gamma$ would be the fraction between the
$\ell_2$-norm of the vector $\text{diag}(\Tmat_{i,j})$ and the
$\ell_2$-norm of the vector $\mvect_{i}\odot\mvect_{j}$, i.e.,
$\frac{\|\text{diag}(\Tmat_{i,j})\|_2}{\|\mvect_{i}\odot\mvect_{j}\|_2}$.
Therefore, $\gamma$ is large if
$\gamma > \mu_{2}(\Tmat(\lambda)) = \max \{\frac{\|\text{diag}(\Tmat_{i,j})\|_2}{\|\mvect_{i}\odot\mvect_{j}\|_2}\}_{1\leq i,j \leq 2}$.
$\gamma$ is small for the opposite case.

In general, at the beginning of the optimization, we can select a small
$\lambda$ and a small $\gamma$ to relax both the PDE constraint and the
rank-1 constraint. As the optimization proceeds, we can increase
$\lambda$ and $\gamma$ to strengthen both constraints so that the
solution can converge to the optimal solution of conventional FWI.

\subsection{Computational cost
analysis}\label{computational-cost-analysis}

The major computational cost of the proposed LRWI is to invert the
$2\ngrid \times 2\ngrid$ matrix $\Smatt^{\top}\Smatt$ in
Equation~\ref{AnaSolution} to obtain $\uvect^{\star}$. If we use a
direct solver to invert $\Smatt^{\top}\Smatt$, the computational cost
will be $\mathcal{O}(8\ngrid^3)$. With $\uvect^{\star}$ in hand, the
computation of the gradients does not include additional matrix
inverses. At each iteration, we alternately update $\mvect$ and
$\theta$. Therefore, the total computational cost for LRWI is
$\mathcal{O}(16\ngrid^3\nfreq)$ for each iteration. Compared to
conventional FWI, whose computational cost is
$\mathcal{O}(2\ngrid^3\nfreq)$ for each iteration, LRWI is eight times
expensive. Considering the increased computational cost, instead of
using LRWI for the whole inversion, we suggest using LRWI to create a
better initial model for FWI.

\section{Numerical examples}\label{numerical-examples}

To investigate the feasibility of the proposed LRWI approach, we conduct
numerical examples on two well-known models i.e.~the Marmousi model
\citep{versteeg1994marmousi} and the Overthrust model
\citep{overthrust1996}. In both examples, we will study the performances
of the proposed LRWI for different selections of $\lambda$ and $\gamma$,
and investigate the performance with respect to the starting frequency.

\subsection{Marmousi model}\label{marmousi-model}

\begin{figure}
\centering
\subfloat[True model
$\mvec_{\text{t}}$\label{fig:Mtrue}]{\includegraphics[width=0.500\hsize]{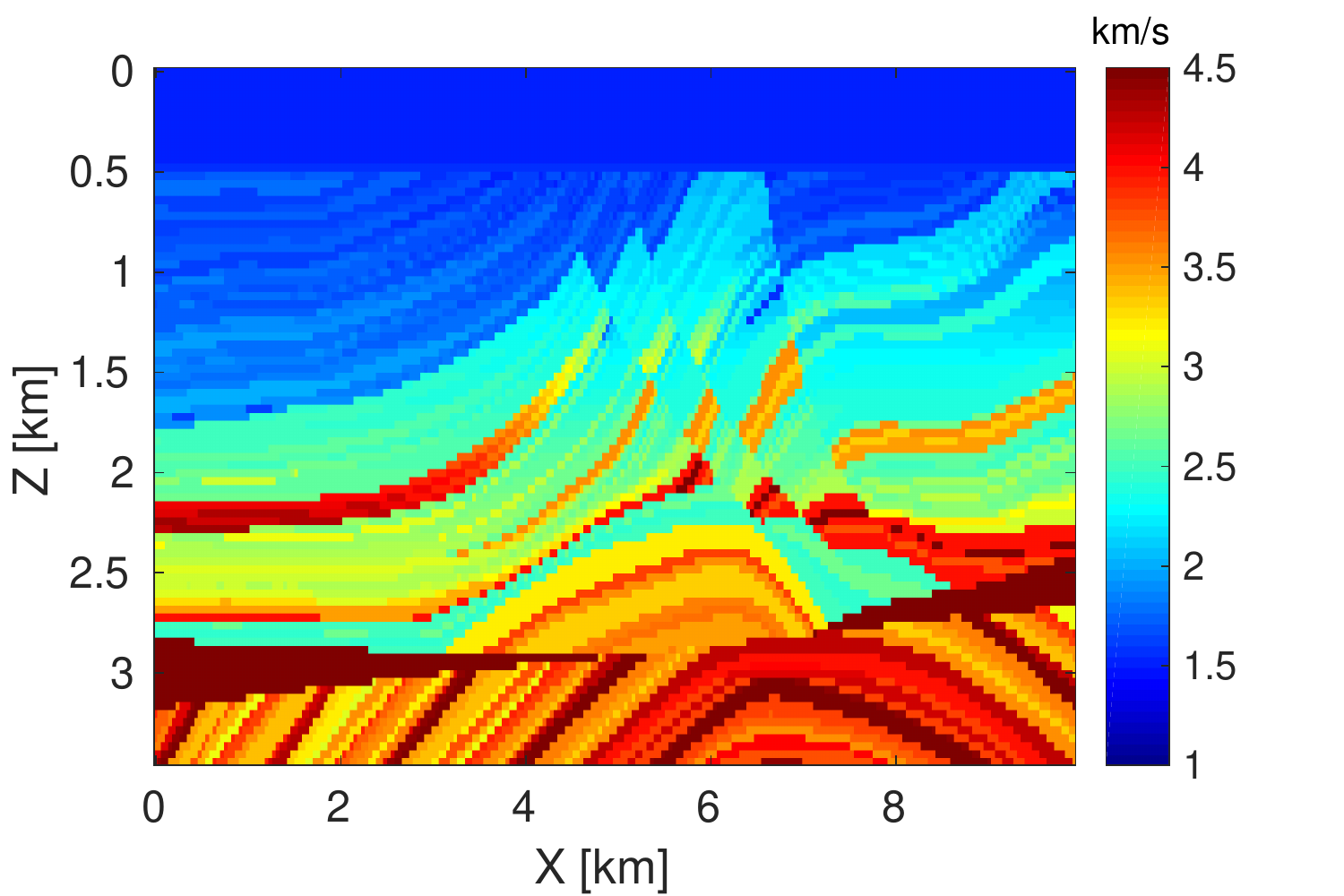}}
\subfloat[Initial model
$\mvec^{(0)}$\label{fig:Mini1}]{\includegraphics[width=0.500\hsize]{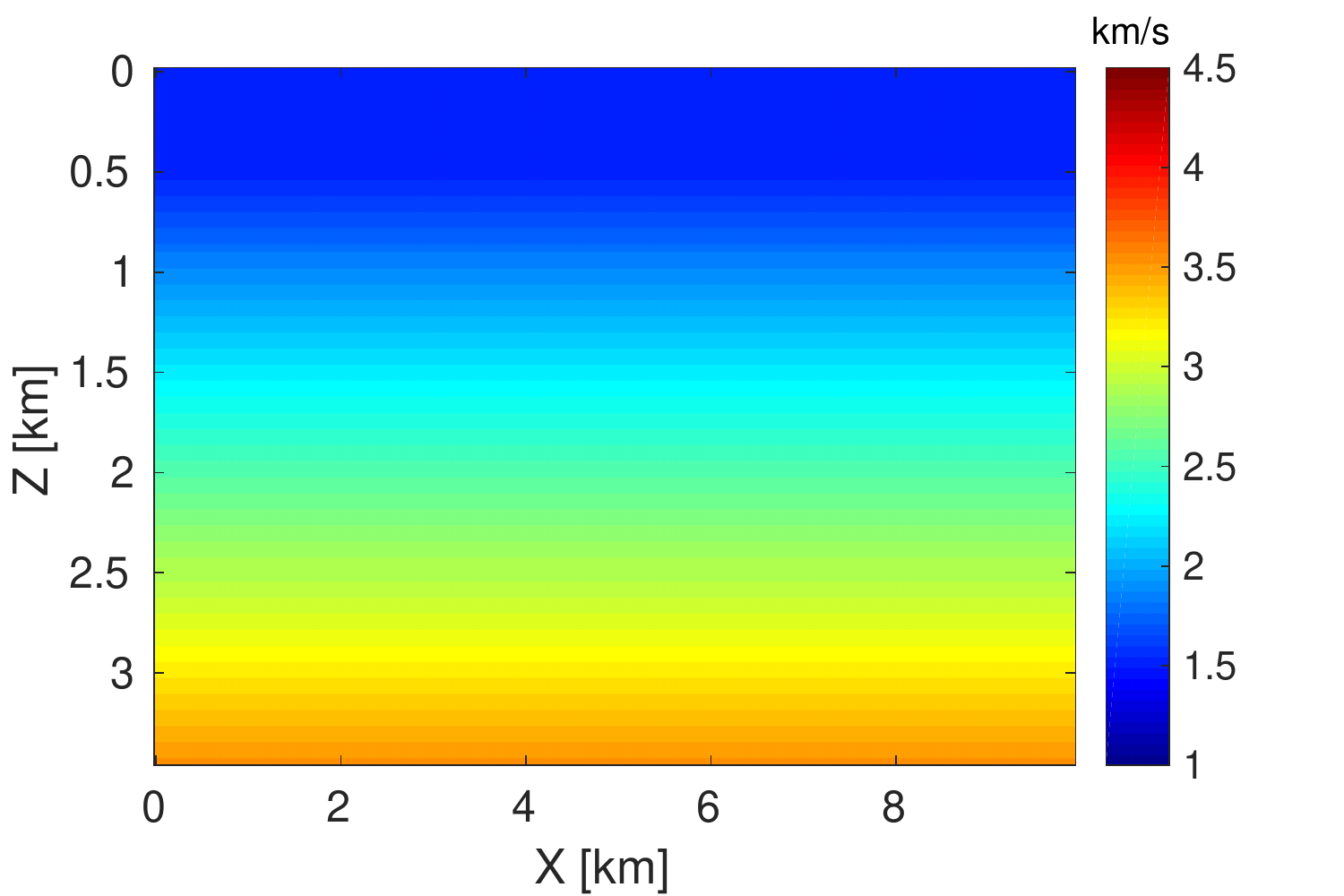}}
\caption{ (a) The true velocity model; (b) The initial velocity
model.}\label{fig:Marmousi}
\end{figure}

We first conduct an example on the Marmousi-2 model $\mvec_\text{t}$
shown in Figure~\ref{fig:Mtrue}. We use a Ricker wavelet centered at
$15 \mHz$ to simulate 49 sources at the depth of $z = 0.04 \mkm$ with a
sampling interval of $0.5 \mkm$. The data are recorded by 247 receivers
at the same depth with a sampling interval of $0.04 \mkm$. As is
commonly practiced, we perform the frequency continuation
\citep{bunks1995multiscale} using three frequency bands of
$\{2.0, 2.5,3.0\}\mHz,\, \{5.0, 6.0, 7.0\}\mHz $, and
$\{7.0, 8.0, 9.0\}\mHz $. We discretize the model with $0.04 \mkm$
grids. We compare the performances of conventional FWI, conventional
WRI, and the proposed LRWI. For conventional FWI and WRI, we use the
l-BFGS method to solve the optimization problem, while we use
Algorithm~\ref{alg:1} to solve the LRWI. Due to the computational cost,
we use LRWI to conduct an inversion on the lowest frequency band and
then use the obtained model as the initial model for conventional FWI.
All three approaches use 45 iterations for each frequency band.

To initialize the inversion, we conduct FWI and WRI with the 1D
monotonously increasing velocity model $\mvec^{(0)}$ shown in
Figure~\ref{fig:Mini1}. For LRWI we select
$\theta^{(0)} = \frac{\pi}{4}$ and
$\mvect^{(0)} = (\sin\theta^{(0)}\mvec^{(0)}, \cos\theta^{(0)}\mvec^{(0)})$.
We conduct conventional WRI with four different selections of the
penalty parameter $\lambda$, i.e.
$\lambda = \beta_{1}\mu_{1}(\Amat^{\top}\Pmat^{\top}\Pmat\Amat^{-1})$
with $\beta_1 = 1e$-8, $1e$-$4, 1e0, \text{and}\, 1e4$. For the proposed
LRWI, we use the same selection for $\lambda$ and select six different
$\gamma$'s for each $\lambda$. We select
$\gamma = \beta_2\mu_2(\Tmat(\lambda))$, with
$\beta_2 = 1e\text{-}16, 1e\text{-}12, 1e\text{-}8, 1e\text{-}4, 1e0, \text{and} 1e4$.
The selections of $\beta_1$ and $\beta_2$ can not be extremely small,
otherwise the matrix
$\Smatt(\beta_1,\beta_2)^{\top}\Smatt(\beta_1,\beta_2)$ would be close
to singular or badly scaled.

\begin{figure}
\centering
\includegraphics[width=0.500\hsize]{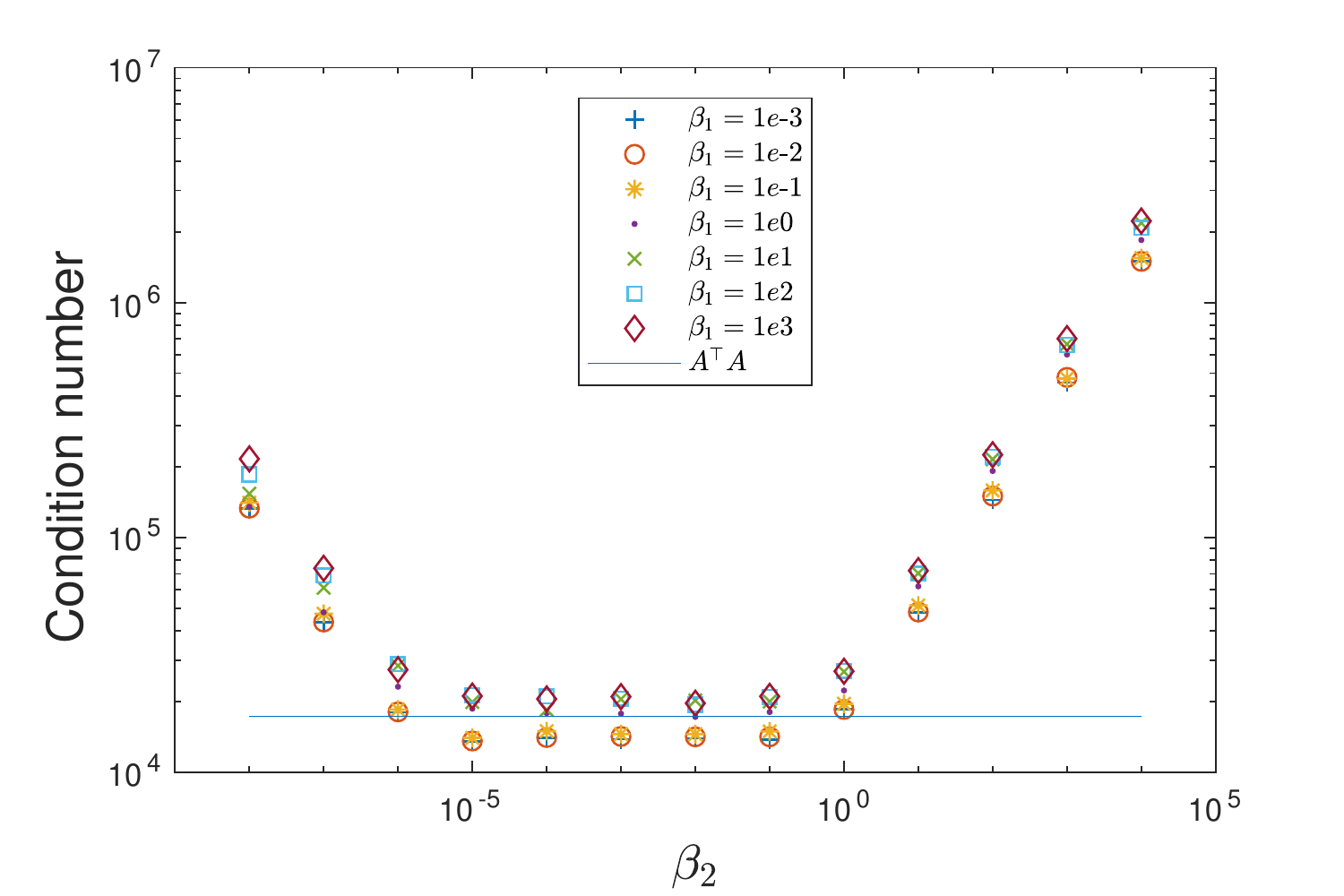}
\caption{Condition number of the matrix $\Smatt^{\top}\Smatt$ versus the
values of $\beta_1$ and $\beta_2$.}\label{fig:MCN}
\end{figure}

Before the inversion, we first study the condition number of the matrix
$\Smatt^{\top}\Smatt$ with respect to the selection of $\beta_1$ and
$\beta_2$. We use the initial model $\mvec^{(0)}$ to form the Helmholtz
matrix $\Amat(\mvec^{(0)})$ and compute the condition number of the
matrix $\Amat(\mvec^{(0)})^{\top}\Amat(\mvec^{(0)})$ as a reference
(c.f. the blue line in Figure~\ref{fig:MCN}). Then we use the initial
model $\mvect^{(0)}$ and different selections of $\beta_1$ and $\beta_2$
to form the matrix
$\Smatt(\beta_1,\beta_2)^{\top}\Smatt(\beta_1,\beta_2)$. The ranges for
$\beta_1$ and $\beta_2$ are {[}$1e$-3, $1e3${]} and {[}$1e$-8, $1e4${]},
respectively. The condition number of the matrix
$\Smatt(\beta_1,\beta_2)^{\top}\Smatt(\beta_1,\beta_2)$ with respect to
different selections of $\beta_1$ and $\beta_2$ are plotted in
Figure~\ref{fig:MCN}. We can observe that when
$1e\text{-}6\leq\beta_2\leq1e0$, the condition number of
$\Smatt(\beta_1,\beta_2)^{\top}\Smatt(\beta_1,\beta_2)$ is close to that
of $\Amat(\mvec^{(0)})^{\top}\Amat(\mvec^{(0)})$. When
$\beta_2 \leq 1e\text{-}6$, the condition number increases 10 times as
$\beta_2$ decreases 100 times. When $\beta_2 \geq 1e0$, the condition
number increases 10 times as $\beta_2$ increases 100 times. Compared to
$\beta_2$, $\beta_1$ possesses a less influence to the condition number
of $\Smatt(\beta_1,\beta_2)^{\top}\Smatt(\beta_1,\beta_2)$. The
variation of the condition number of
$\Smatt(\beta_1,\beta_2)^{\top}\Smatt(\beta_1,\beta_2)$ with respect to
$\beta_1$ is less than that of $\beta_2$.

\begin{figure}
\centering
\includegraphics[width=0.500\hsize]{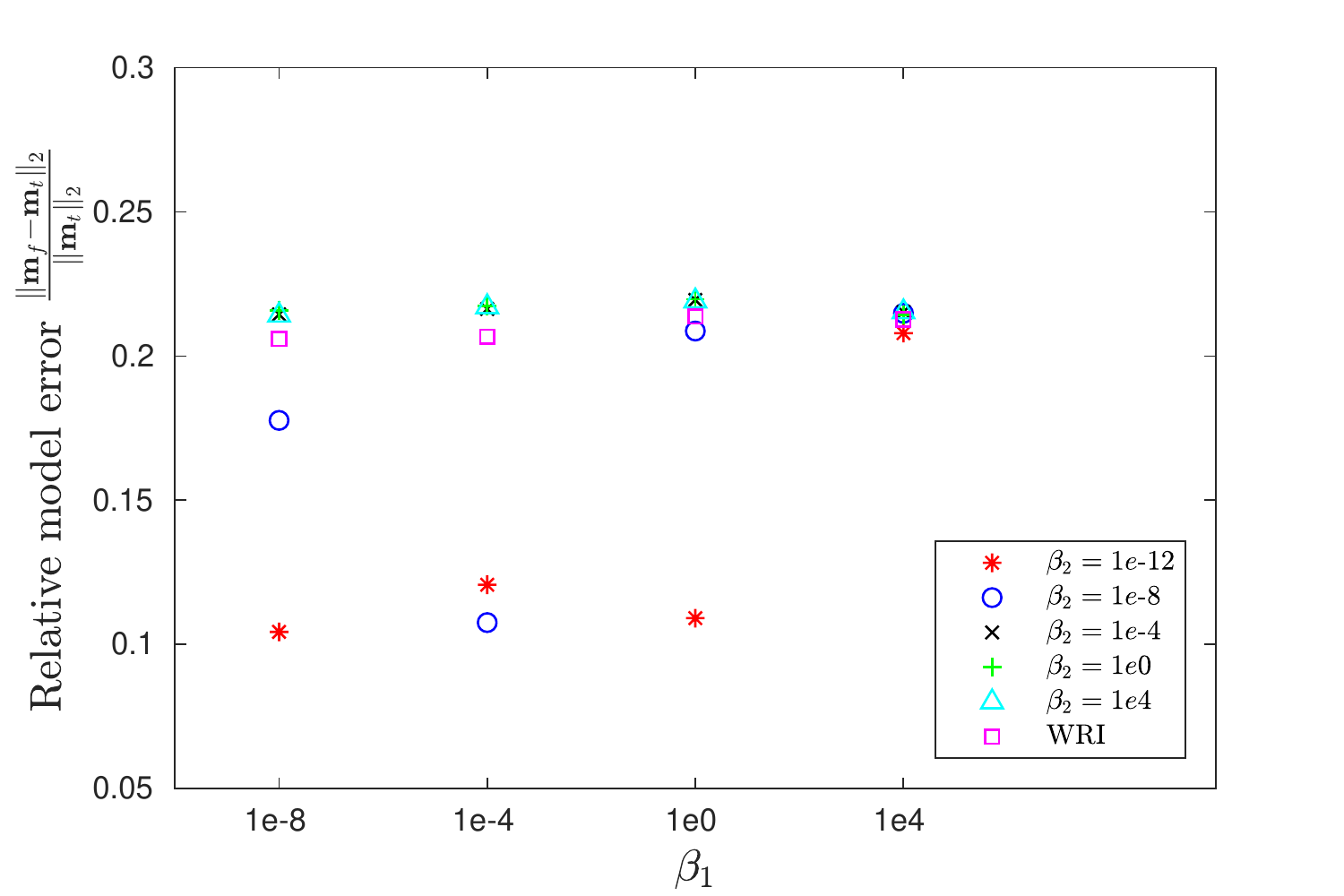}
\caption{Relative model error comparison for WRI with different
selections of $\beta_1$ and LRWI with different selections of $\beta_1$
and $\beta_2$.}\label{fig:MMEpen}
\end{figure}

\begin{figure}
\centering
\subfloat[Result of FWI after the first frequency
band\label{fig:MRfwi01}]{\includegraphics[width=0.300\hsize]{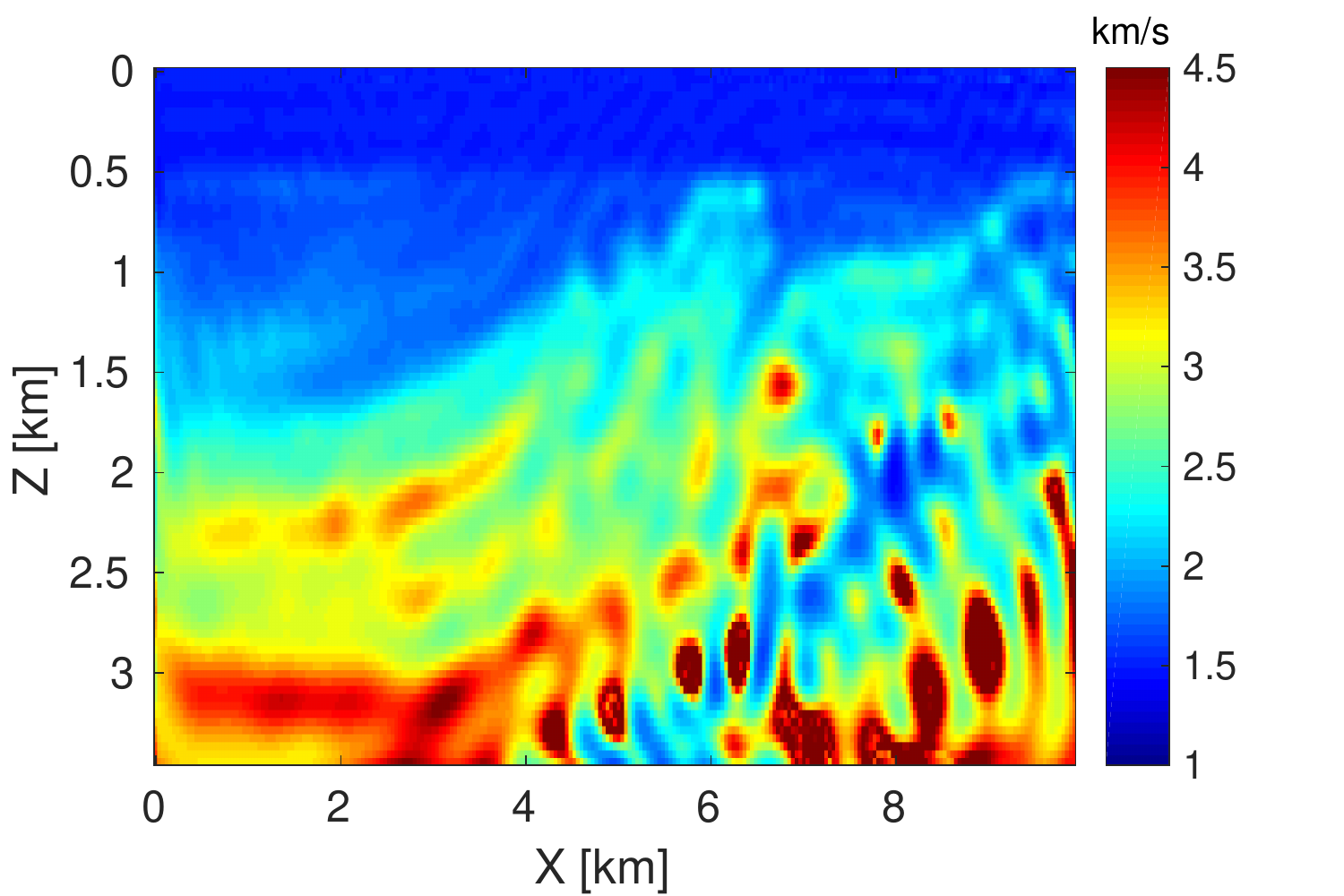}}
\subfloat[Result of WRI after the first frequency
band\label{fig:MRwri01}]{\includegraphics[width=0.300\hsize]{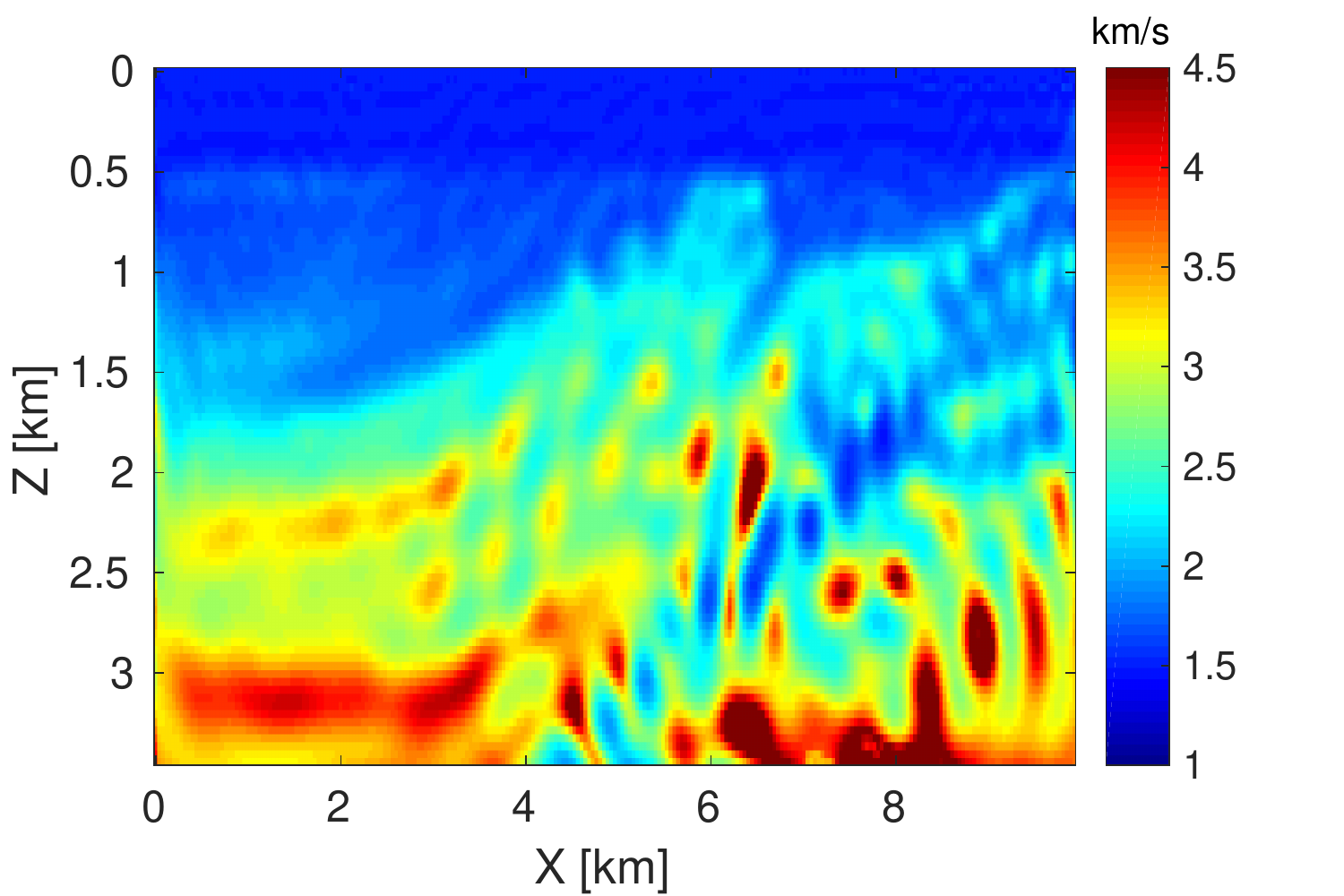}}
\subfloat[Result of LRWI after the first frequency
band\label{fig:MRlrwi01}]{\includegraphics[width=0.300\hsize]{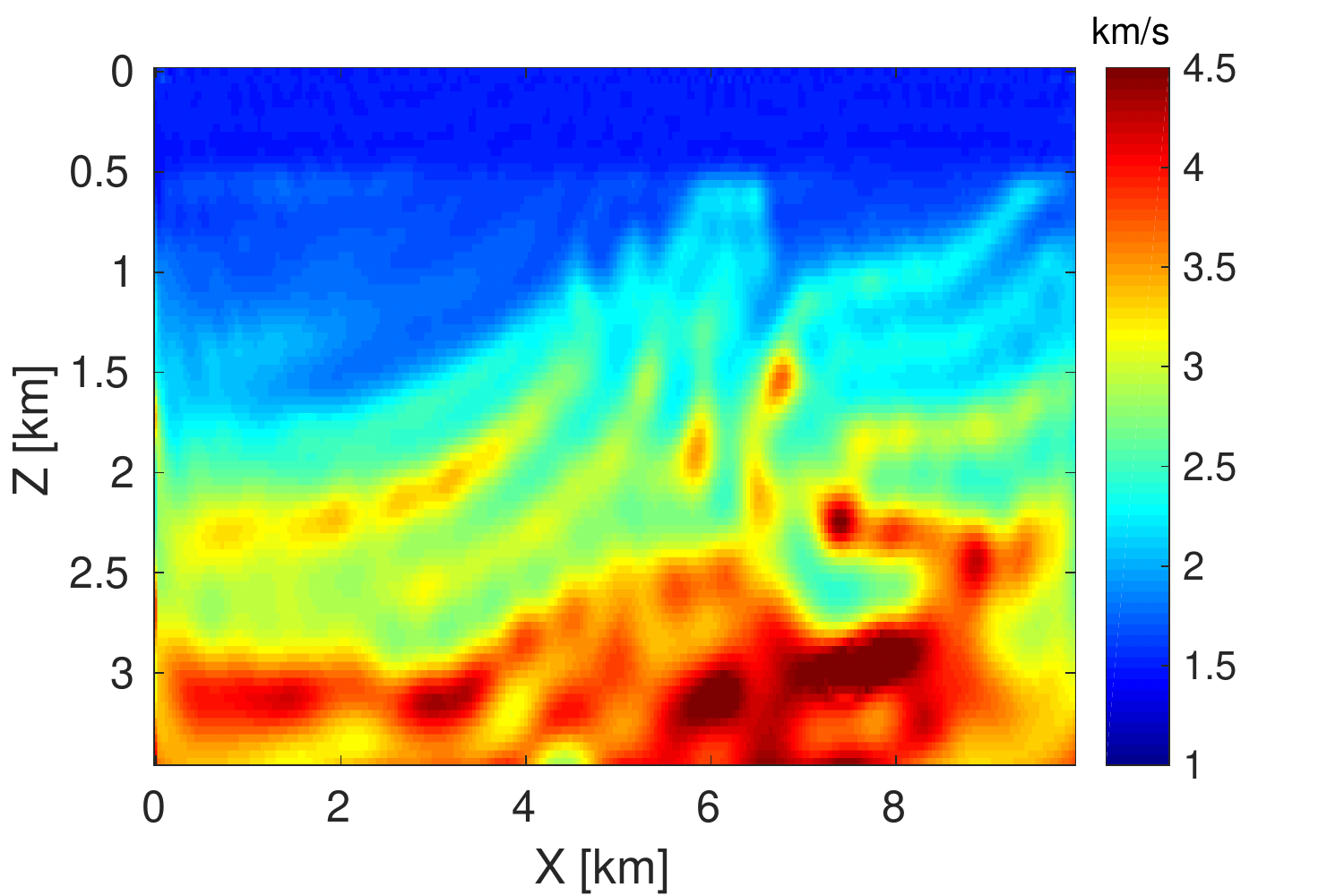}}
\\
\subfloat[Final result of
FWI\label{fig:MRfwi03}]{\includegraphics[width=0.300\hsize]{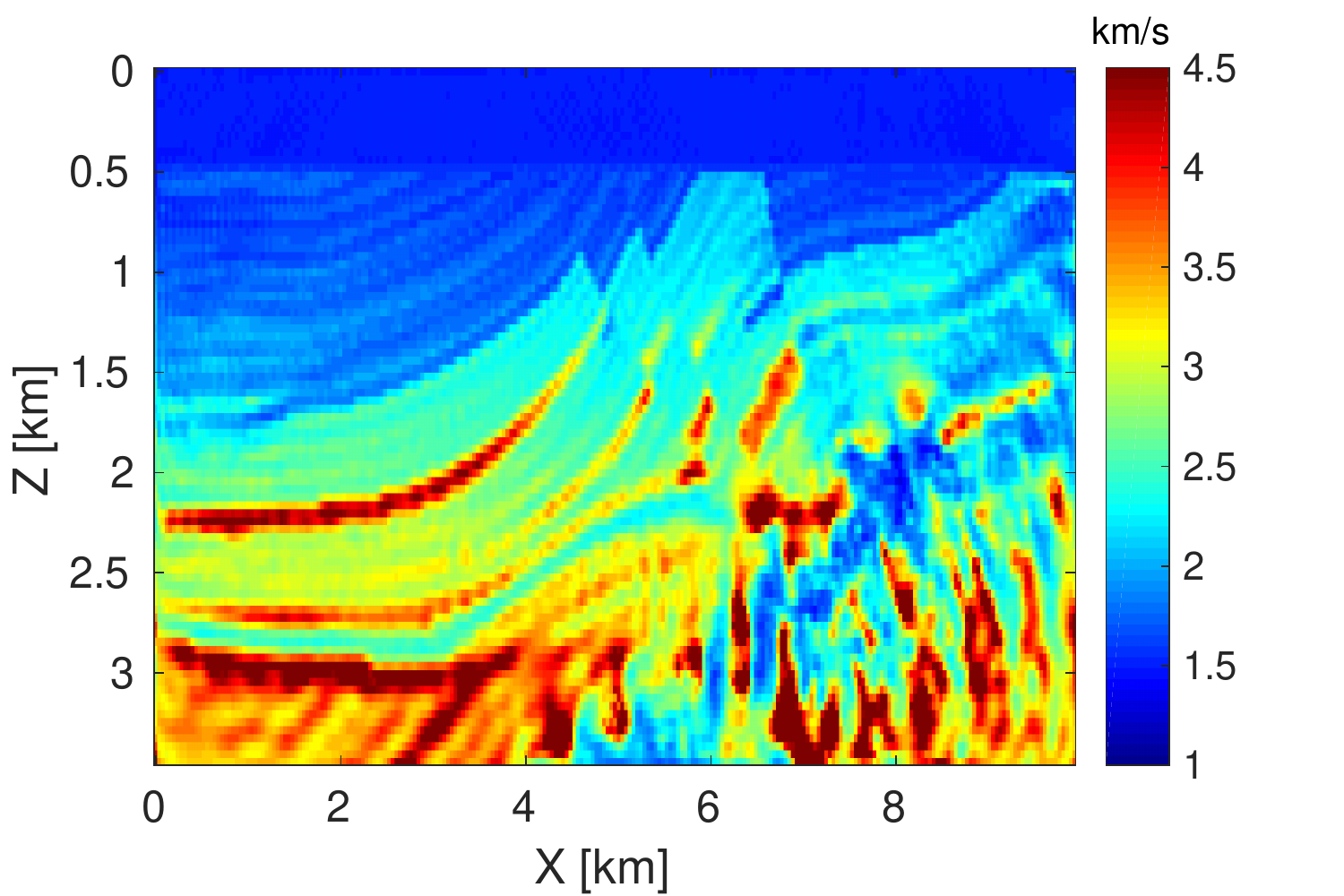}}
\subfloat[Final result of
WRI\label{fig:MRwri03}]{\includegraphics[width=0.300\hsize]{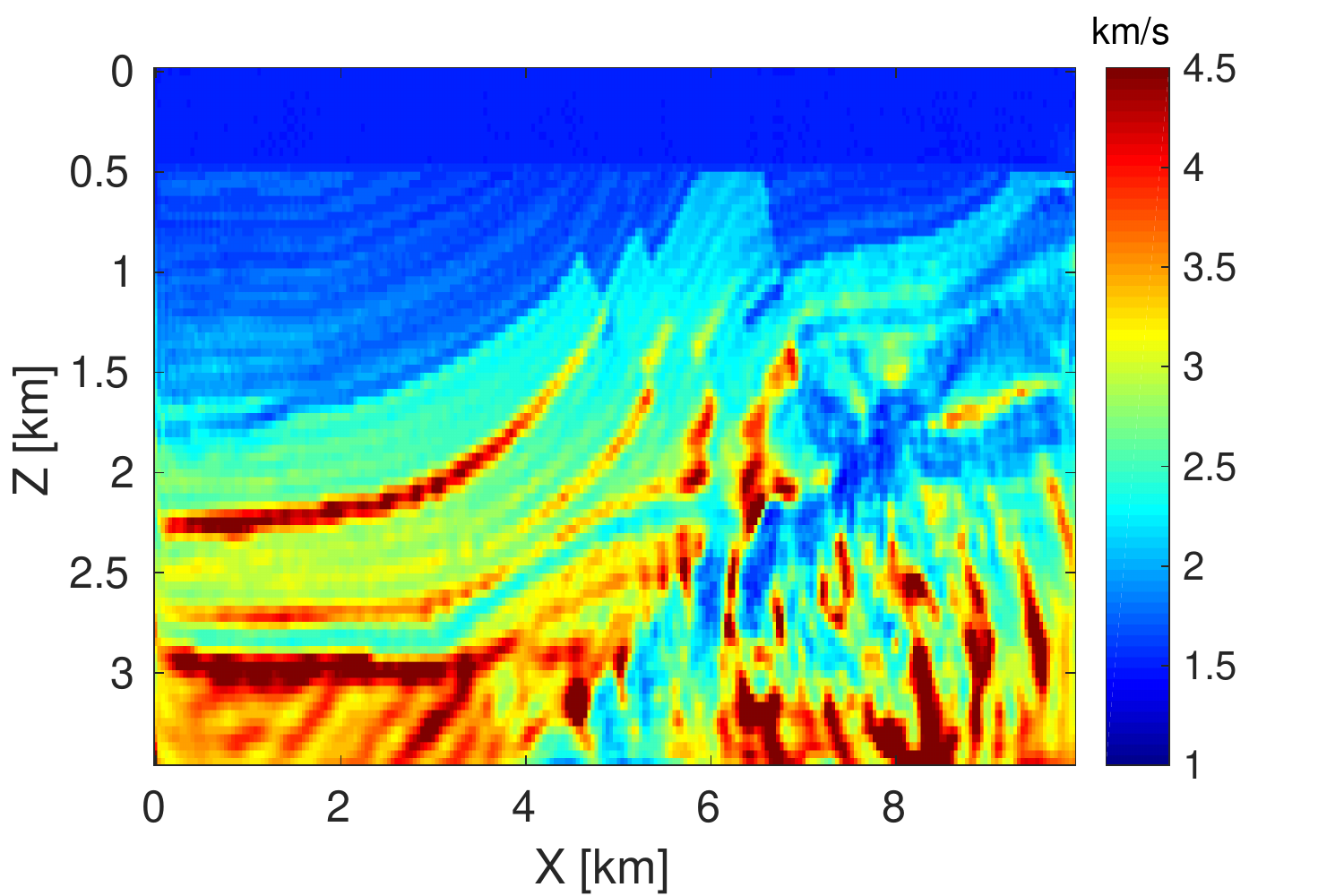}}
\subfloat[Final result of
LRWI\label{fig:MRlrwi03}]{\includegraphics[width=0.300\hsize]{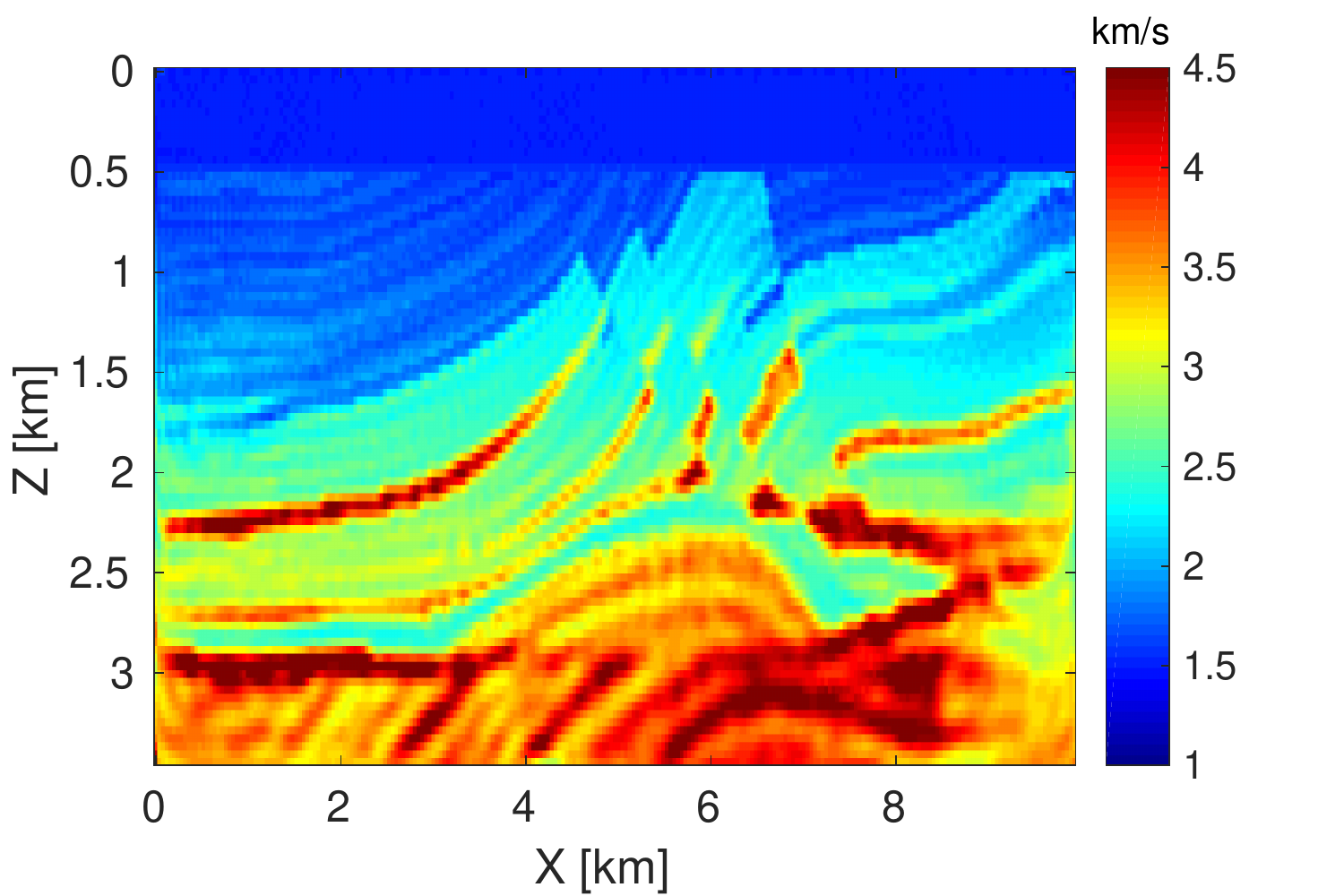}}
\caption{(a) - (c) Results of FWI, WRI, and LRWI after the first
frequency band. (d) - (f) Final results of FWI, WRI, and
LRWI.}\label{fig:MResult0}
\end{figure}

Figure~\ref{fig:MMEpen} shows the relative model error
$\frac{\|\mvec_\text{t}-\mvec_\text{f}\|_2}{\|\mvec_\text{t}\|_2}$
between the true model $\mvec_\text{t}$ and the final inverted model
$\mvec_\text{f}$ obtained by conventional WRI and LRWI with different
selections of $\beta_1$ and $\beta_2$. According to
Figure~\ref{fig:MMEpen}, the selection of $\beta_1=1e$-4 produces the
best result for conventional WRI, and the selection of $\beta_1=1e$-8
and $\beta_2=1e$-12 produces the best result for the LRWI.
Figure~\ref{fig:MResult0} shows the final inverted models
$\mvec_\text{f}$ of conventional FWI, WRI with the best selection of
$\beta_1$, and LRWI with the best selection of $\beta_1$ and $\beta_2$.
Figures~\ref{fig:MRfwi01} to~\ref{fig:MRlrwi01} show the results of the
three approaches using the data of the first frequency band. Clearly,
under the current experimental settings, both FWI and WRI already
converge to local minima at the first frequency band, despite the fact
that WRI can outperform FWI in some other settings. On the other hand,
LRWI provides a much better model for the following inversion, which
yields a significantly better final result shown in
Figure~\ref{fig:MRlrwi03} compared to those obtained by FWI and WRI
(c.f. Figures~\ref{fig:MRfwi03} and~\ref{fig:MRwri03}).

Figure~\ref{fig:MMEpen} does not include the result of LRWI with
$\beta_2 = 1e$-$16$ due to the fact that the matrix
$\Smatt(\beta_1,\beta_2)^{\top}\Smatt(\beta_1,\beta_2)$ is close to
singular or badly scaled with such a small selection of $\beta_2$.
Therefore, we should avoid selecting too small $\beta_2$ when using
LRWI. Figure~\ref{fig:MMEpen} shows that with the selection of
$\beta_1 \leq 1e0$ and $\beta_2=1e$-$12$ or $\beta_1=1e$-$4$ and
$\beta_2=1e$-$8$, LRWI can reconstruct an inverted model with a relative
model error of $10\%$, which is significantly smaller than those of
conventional WRI and LRWI with other selections of $\beta_1$ and
$\beta_2$. This result implies that when the initial model is poor, LRWI
can bypass the local minima of conventional FWI and WRI by properly
relaxing the wave-equation constraint and the rank-1 constraint.

\begin{figure}
\centering
\includegraphics[width=0.500\hsize]{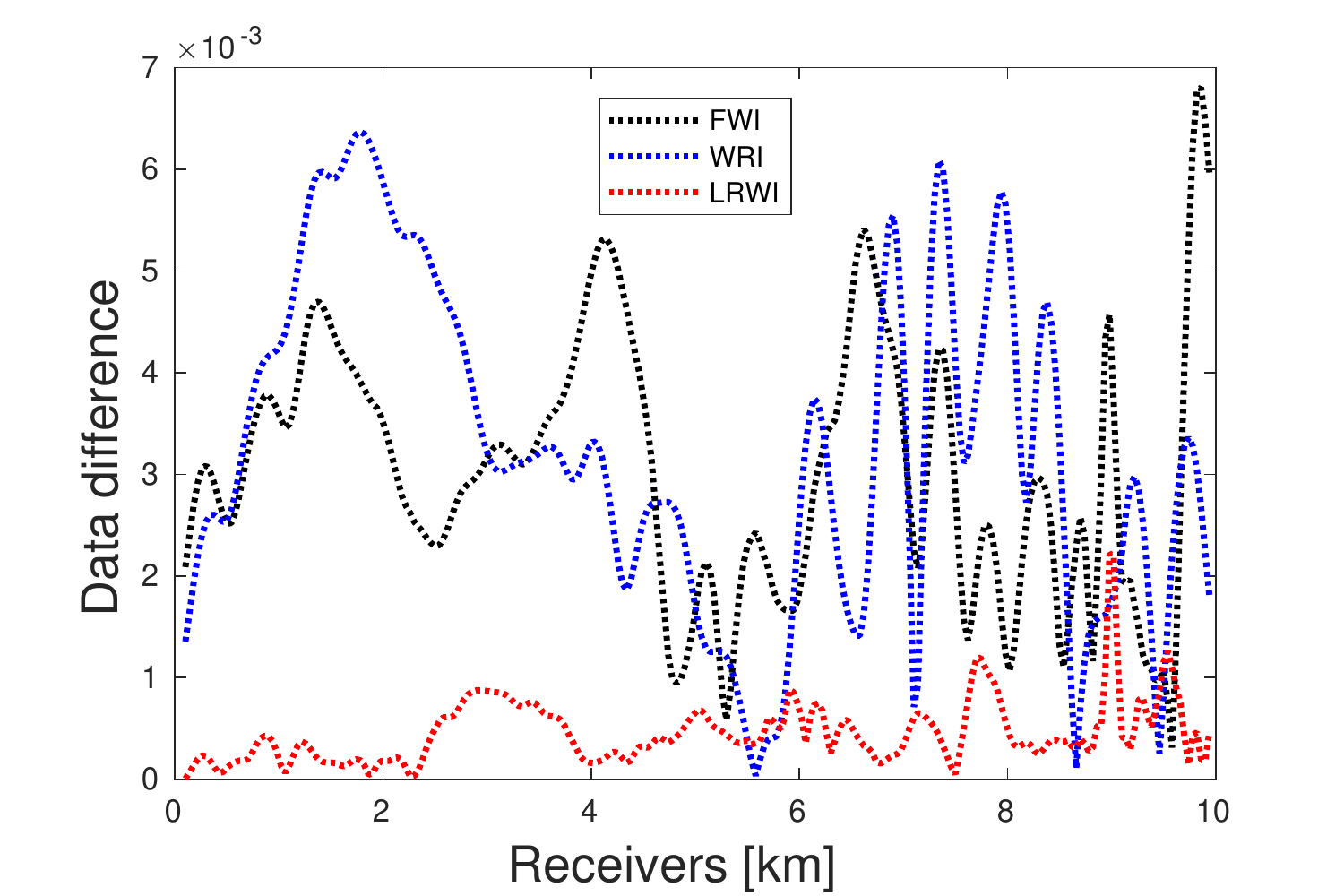}
\caption{Comparison of the absolute data difference
$|\dvec_{\text{obs}} - \dvec_{\text{pred}}|$ for the source located at
$x\, =\, 9\, \mathrm{km}$ and frequency of $3\, \mathrm{Hz}$. The three
lines denote the absolute data differences corresponding to the inverted
results of FWI (black), WRI (blue), and LRWI (red) using the data of the
the first frequency band.}\label{fig:MDEcmp}
\end{figure}

To further compare the inverted results of the three approaches, we use
the three inverted models shown in Figures~\ref{fig:MRfwi01}
-~\ref{fig:MRlrwi01} to compute the predicted data $\dvec_{\text{pred}}$
at the frequency of $3\, \mathrm{Hz}$ for the source located at
$x\, =\, 9\, \mathrm{km}$. We compute the absolute data differences
$|\dvec_{\text{obs}} - \dvec_{\text{pred}}|$ between the observed data
$\dvec_{\text{obs}}$ and the predicted data $\dvec_{\text{pred}}$ and
depict them in Figure~\ref{fig:MDEcmp}. The absolute data difference of
FWI and WRI is more than 6 times larger than that of LRWI. This result
coincides with the fact that FWI and WRI converge to local minima, while
LRWI bypasses the local minima.

\emph{Robustness with respect to the starting frequency} To investigate
the robustness of the three methods with respect to the starting
frequency, we conduct an additional experiment, in which we vary the
starting frequency from 0.5 Hz to 3.0 Hz. We use the same initial model
shown in Figure~\ref{fig:Mini1}. Figure~\ref{fig:OMFreq} illustrates the
relative model errors versus the starting frequency for all the three
methods. According to the previous example, when the relative model
error reaches around 18\%, the inversion converges to a local minimum.
The highest starting frequencies for conventional FWI, WRI, and LRWI to
obtain an inverted model with an acceptable relative model error ($\leq$
14\%) are 1 Hz, 1 Hz, and 2 Hz, respectively. This comparison implies
that under the aforementioned experimental settings, LRWI can conduct a
successful inversion with a starting frequency twice large as that of
conventional FWI and WRI.

\begin{figure}
\centering
\includegraphics[width=0.500\hsize]{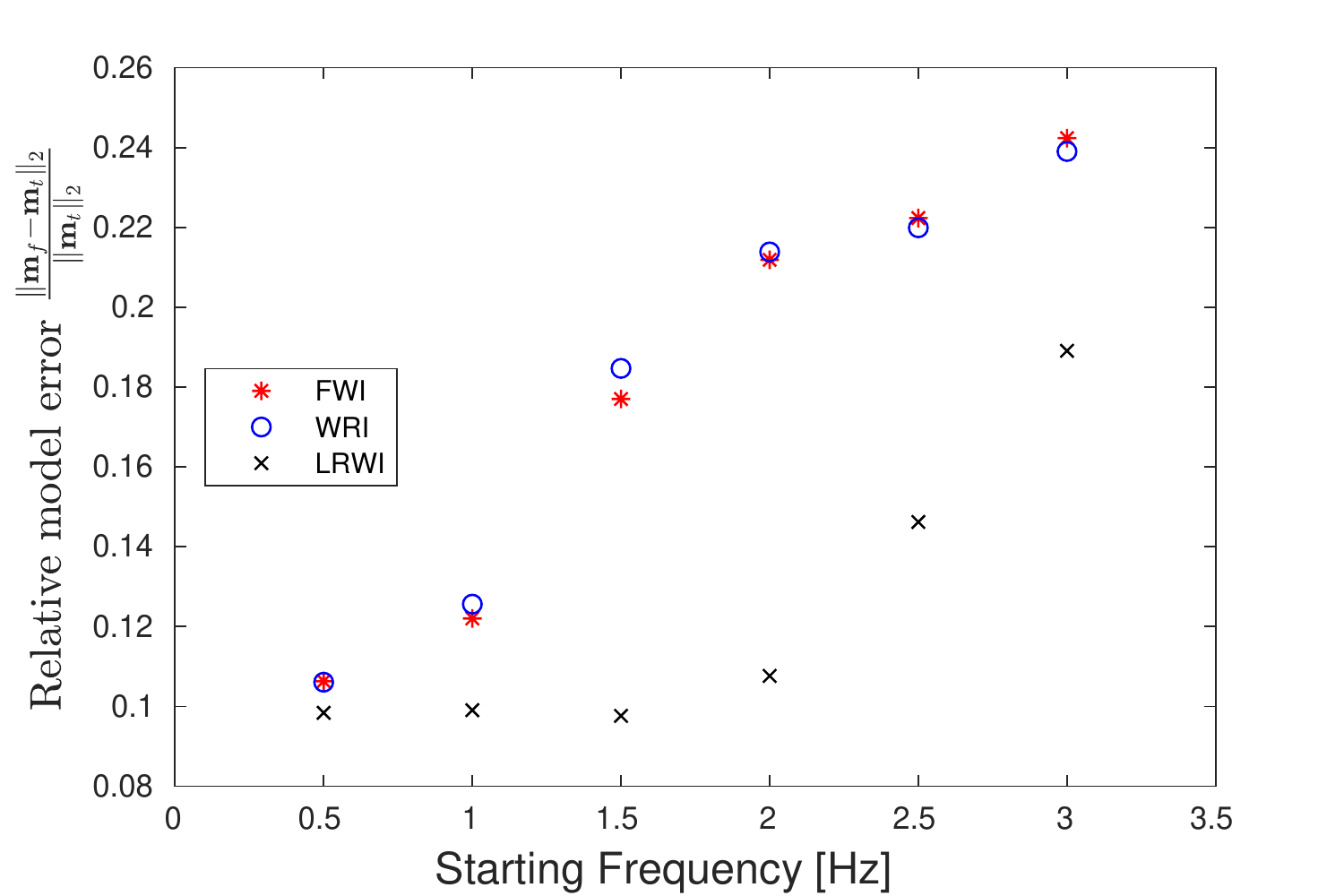}
\caption{Relative model error comparison for FWI($*$), WRI($\circ$), and
LRWI($\times$) using data with different starting
frequencies.}\label{fig:MMFreq}
\end{figure}

\subsection{Overthrust model}\label{overthrust-model}

\begin{figure}
\centering
\subfloat[True model
$\mvec_{\text{t}}$\label{fig:Otrue}]{\includegraphics[width=0.500\hsize]{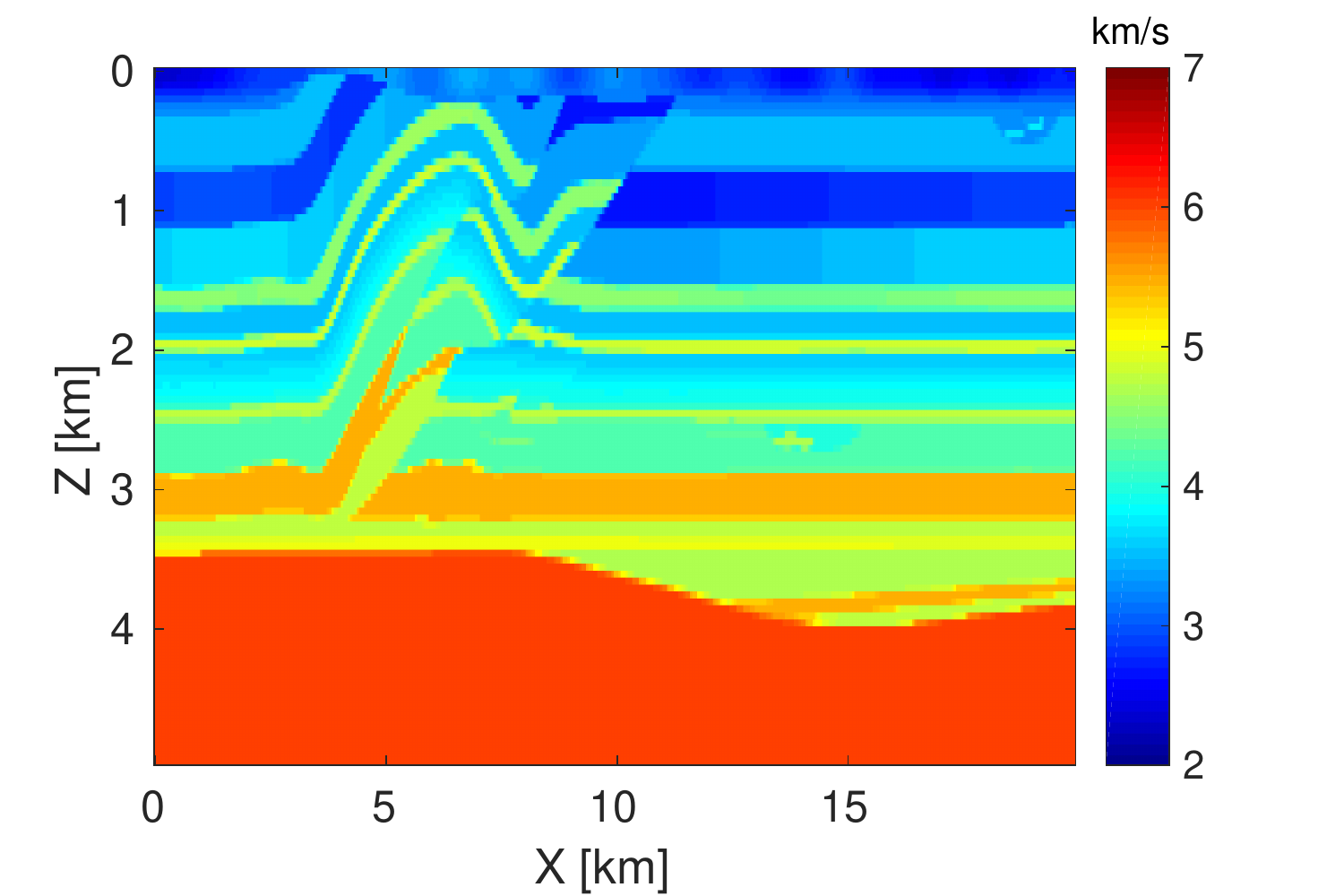}}
\subfloat[Initial model
$\mvec^{(0)}$\label{fig:Oini1}]{\includegraphics[width=0.500\hsize]{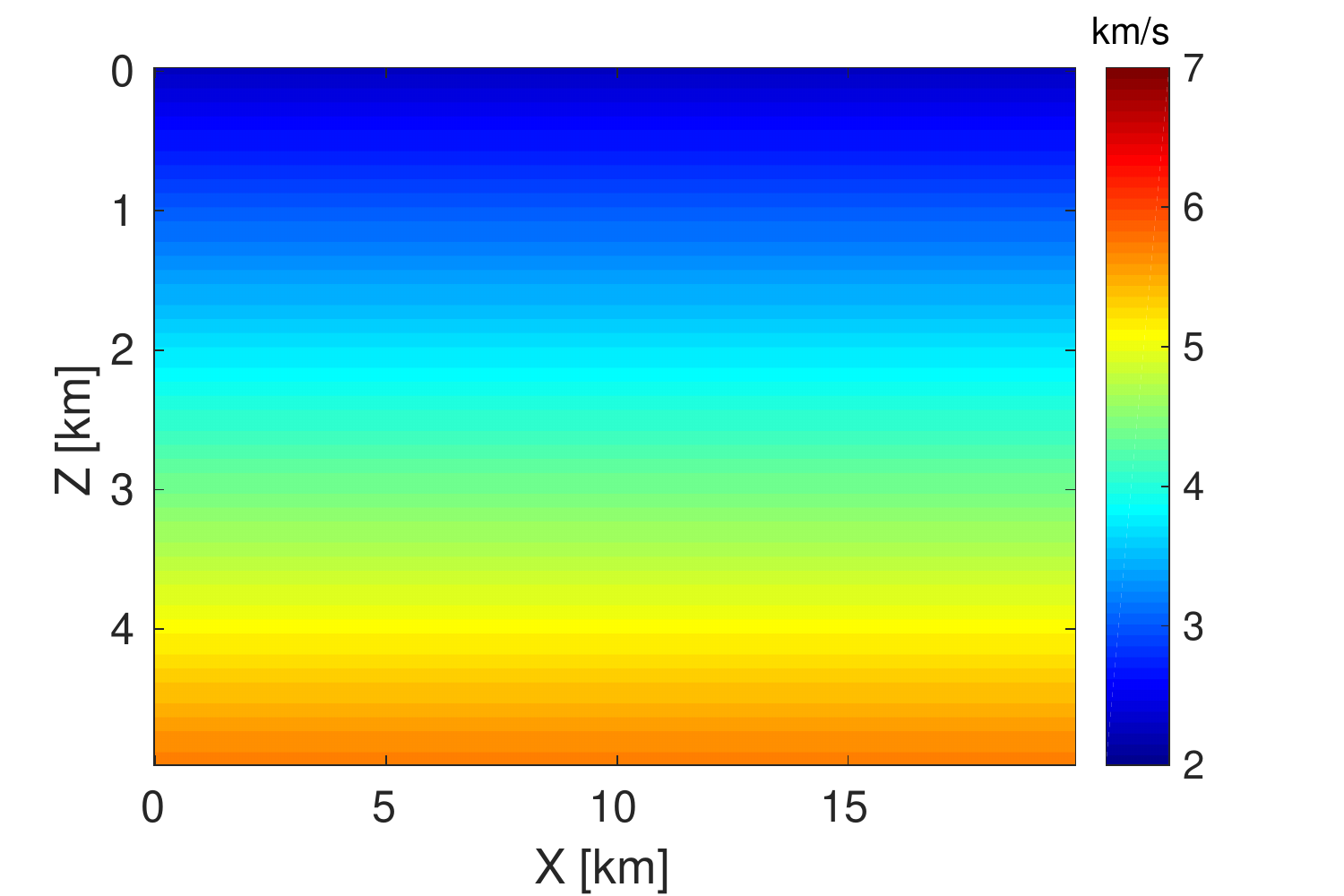}}
\caption{ (a) The true velocity model; (b) The initial velocity
model.}\label{fig:Overthrust}
\end{figure}

We conduct an experiment with the Overthrust model to investigate the
generality of the proposed LRWI with respect to different velocity
structures. Figure~\ref{fig:Otrue} shows the $5\km \times 20\km$
Overthrust model. We place 99 sources and 100 receivers at the depth of
$0.1\km$ with horizontal sampling intervals of $0.2 \km$ and $0.2 \km$,
respectively. As used in the example of the Marmousi model, we conduct
the inversion with the frequency continuation strategy using three
frequency bands of $\{2.0, 2.5,3.0\}\mHz,\, \{5.0, 6.0, 7.0\}\mHz $, and
$\{7.0, 8.0, 9.0\}\mHz $. We discretize the model with $0.05 \km$ grids.
We use the same optimization strategy as that used in the Marmousi
example for the inversion of conventional FWI, conventional WRI, and
LRWI.

We conduct FWI and WRI with the initial model shown in
Figure~\ref{fig:Oini1}. Similar to the previous example, we select
$\theta^{(0)} = \frac{\pi}{4}$ and
$\mvect^{(0)} = (\sin\theta^{(0)}\mvec^{(0)}, \cos\theta^{(0)}\mvec^{(0)})$
to initialize LRWI. We conduct conventional WRI with four different
selections of the penalty parameter $\lambda$, i.e.
$\lambda = \beta_{1}\mu_{1}(\Amat^{\top}\Pmat^{\top}\Pmat\Amat^{-1})$
with $\beta_1 = 1e$-$8,\, 1e$-$4, 1e0, \text{and}\, 1e4$, and select the
one that produces the result with the minimal relative model error as
the output of WRI. For LRWI, we use the same selection for $\lambda$ and
select five different $\gamma$'s for each $\lambda$. We select
$\gamma = \beta_2\mu_2(\Tmat(\lambda))$, with $\beta_2 = 1e$-16,
$1e$-12, $1e$-6, $1e$0, and $1e$6. The combination of $\beta_1$ and
$\beta_2$ that produces the result with the minimal relative model error
is selected as the output of LRWI.

\begin{figure}
\centering
\includegraphics[width=0.500\hsize]{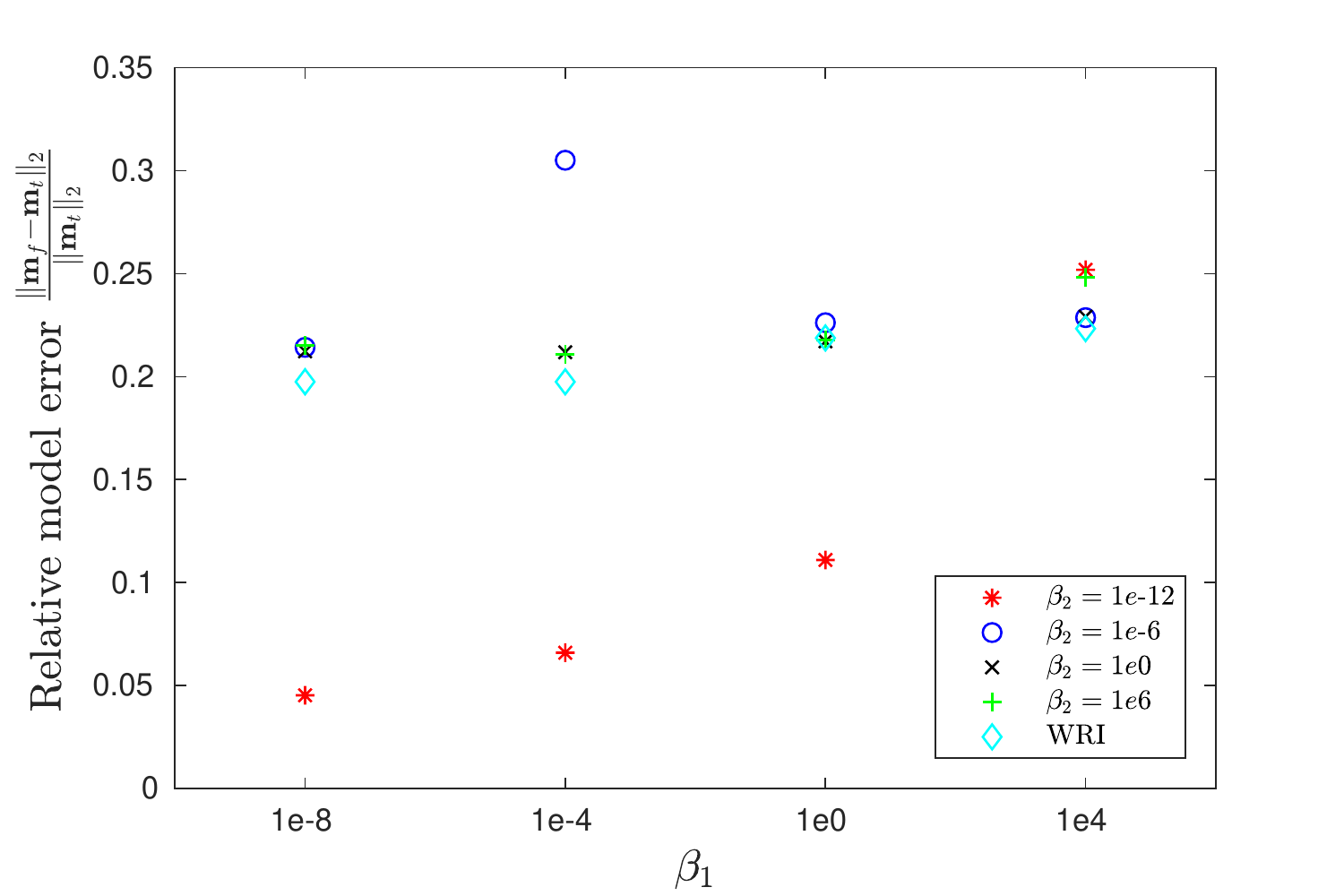}
\caption{Final model error comparison for WRI with different selections
of $\beta_1$ and rank-2 WRI with different selections of $\beta_1$ and
$\beta_2$.}\label{fig:OMEpen}
\end{figure}

\begin{figure}
\centering
\subfloat[Result of FWI after the first frequency
band\label{fig:ORfwi01}]{\includegraphics[width=0.300\hsize]{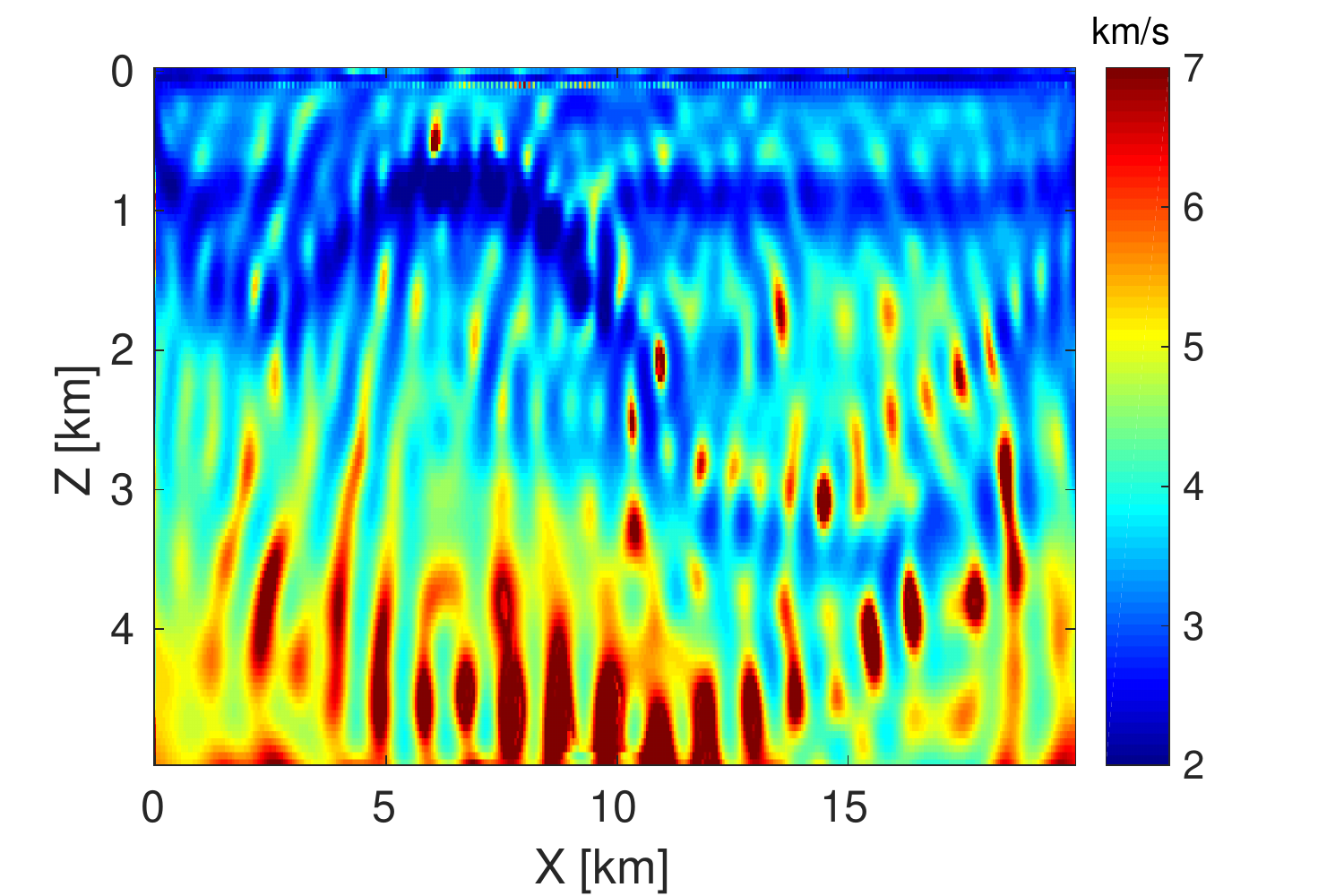}}
\subfloat[Result of WRI after the first frequency
band\label{fig:ORwri01}]{\includegraphics[width=0.300\hsize]{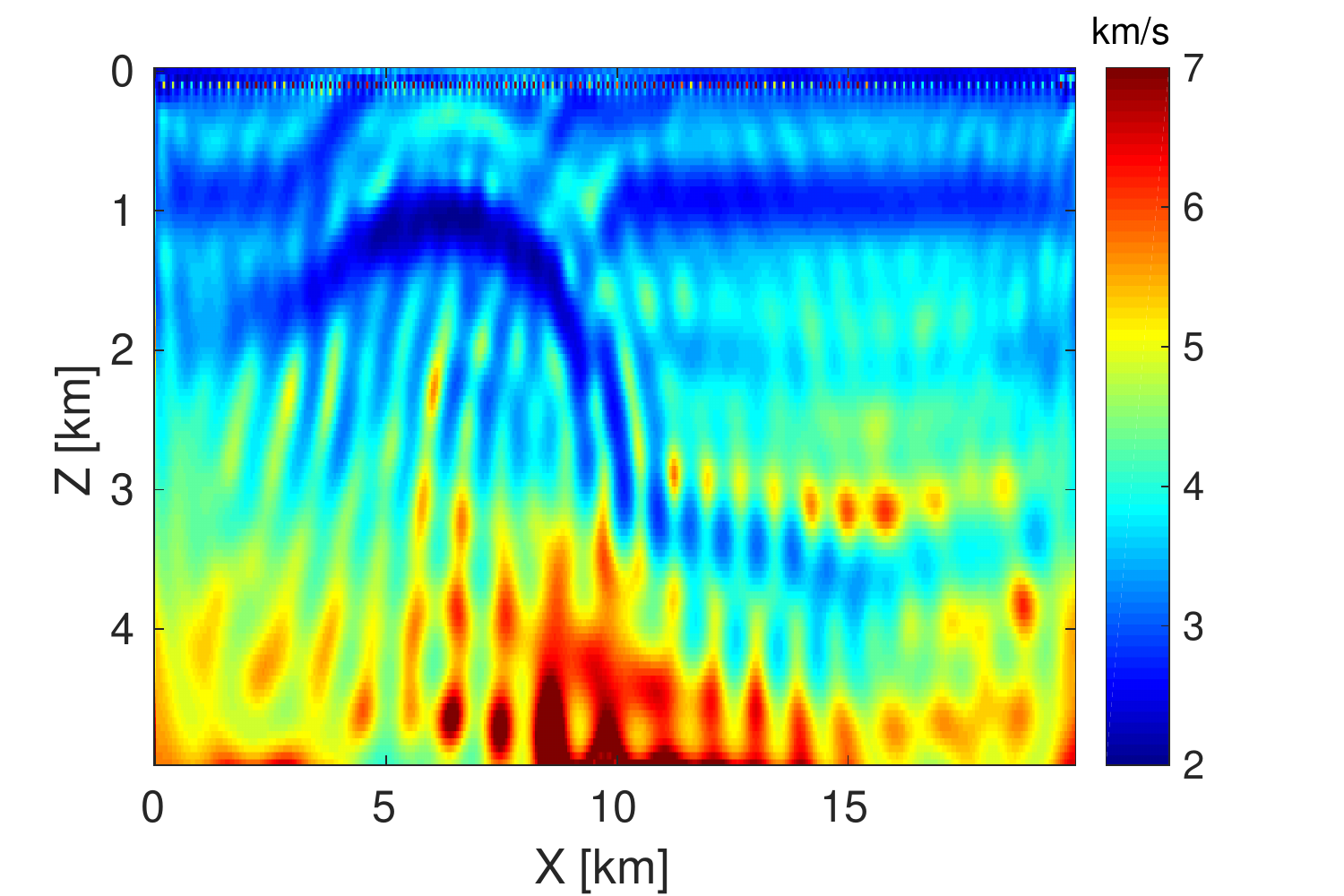}}
\subfloat[Result of LRWI after the first frequency
band\label{fig:ORlrwi01}]{\includegraphics[width=0.300\hsize]{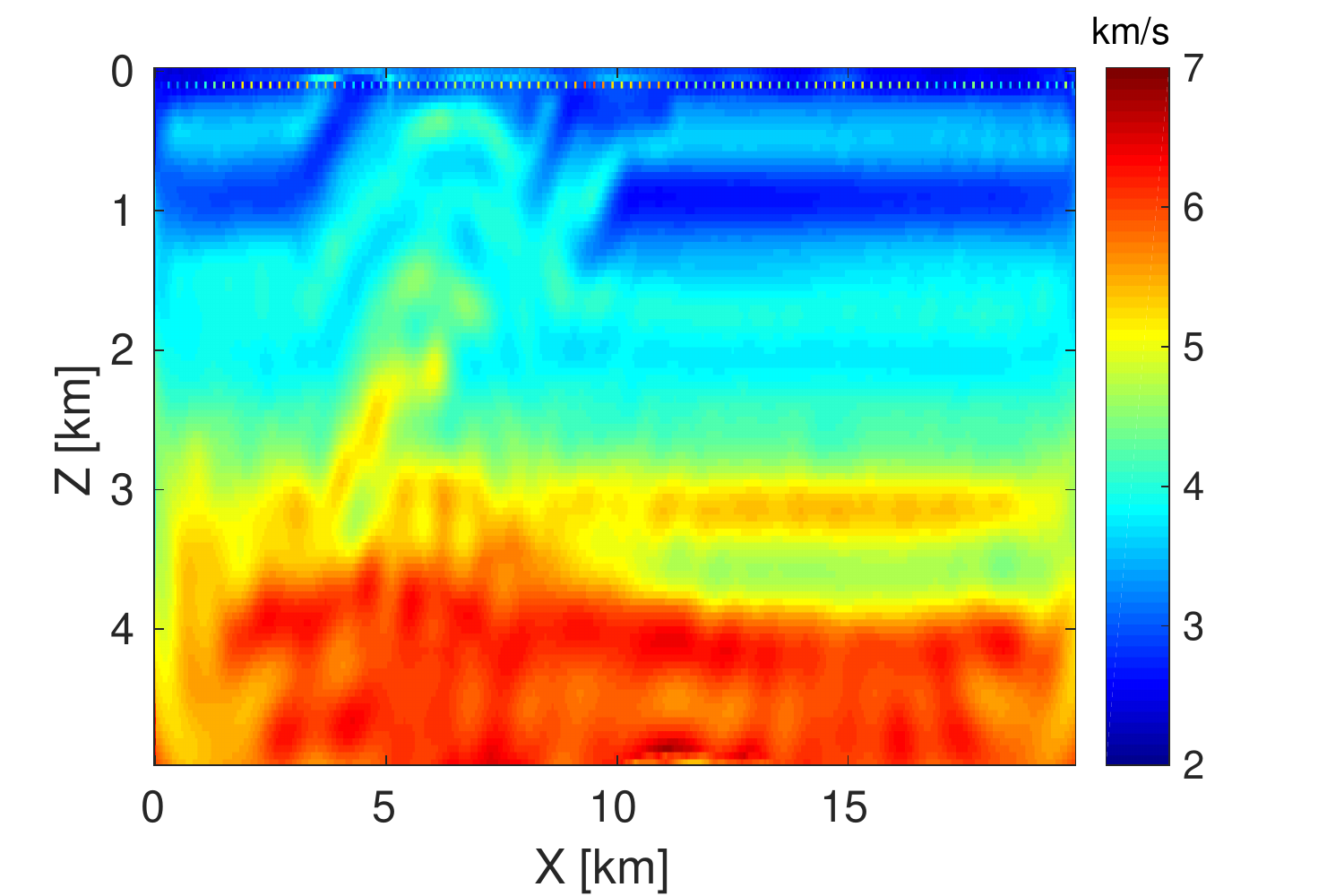}}
\\
\subfloat[Final result of
FWI\label{fig:ORfwi03}]{\includegraphics[width=0.300\hsize]{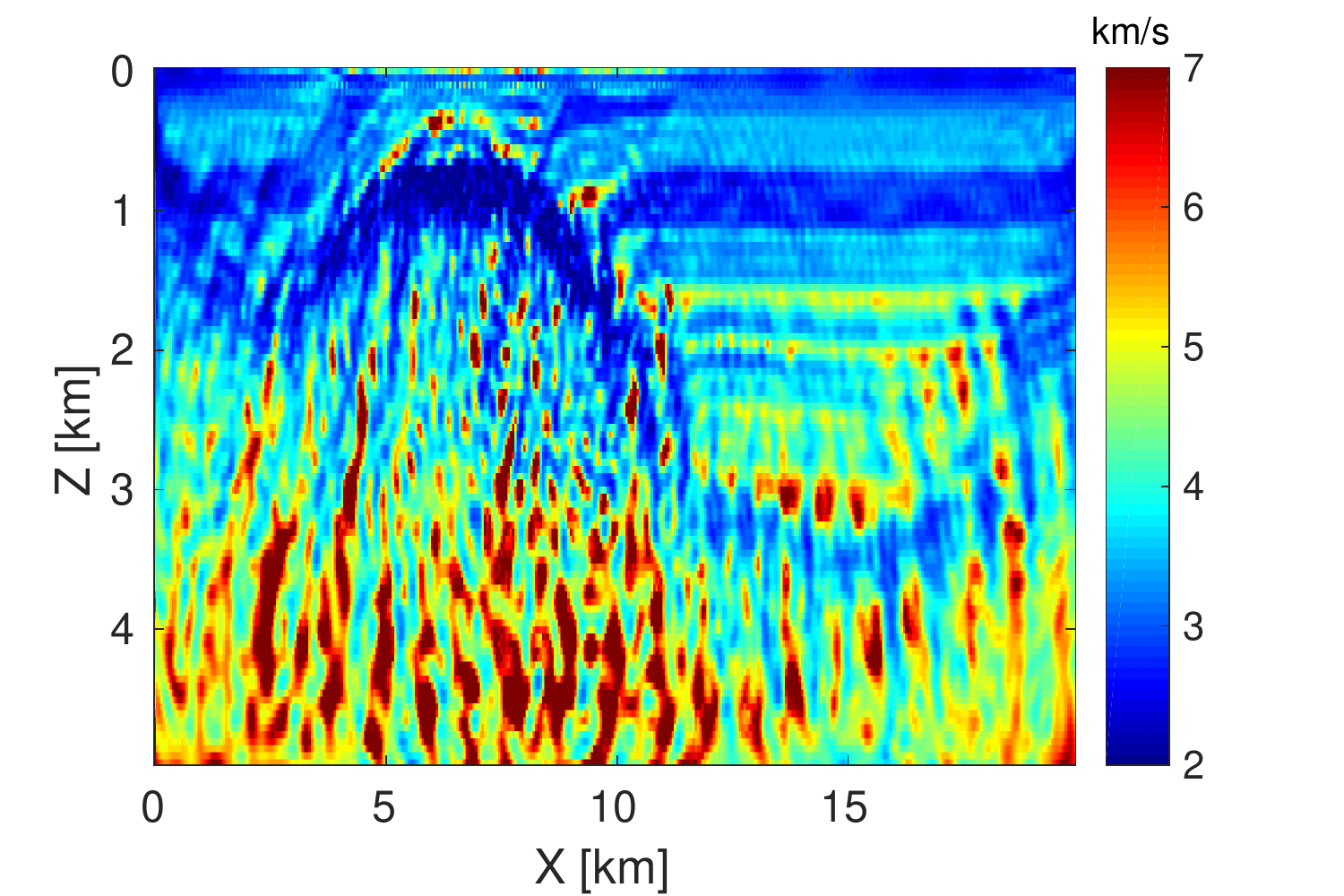}}
\subfloat[Final result of
WRI\label{fig:ORwri03}]{\includegraphics[width=0.300\hsize]{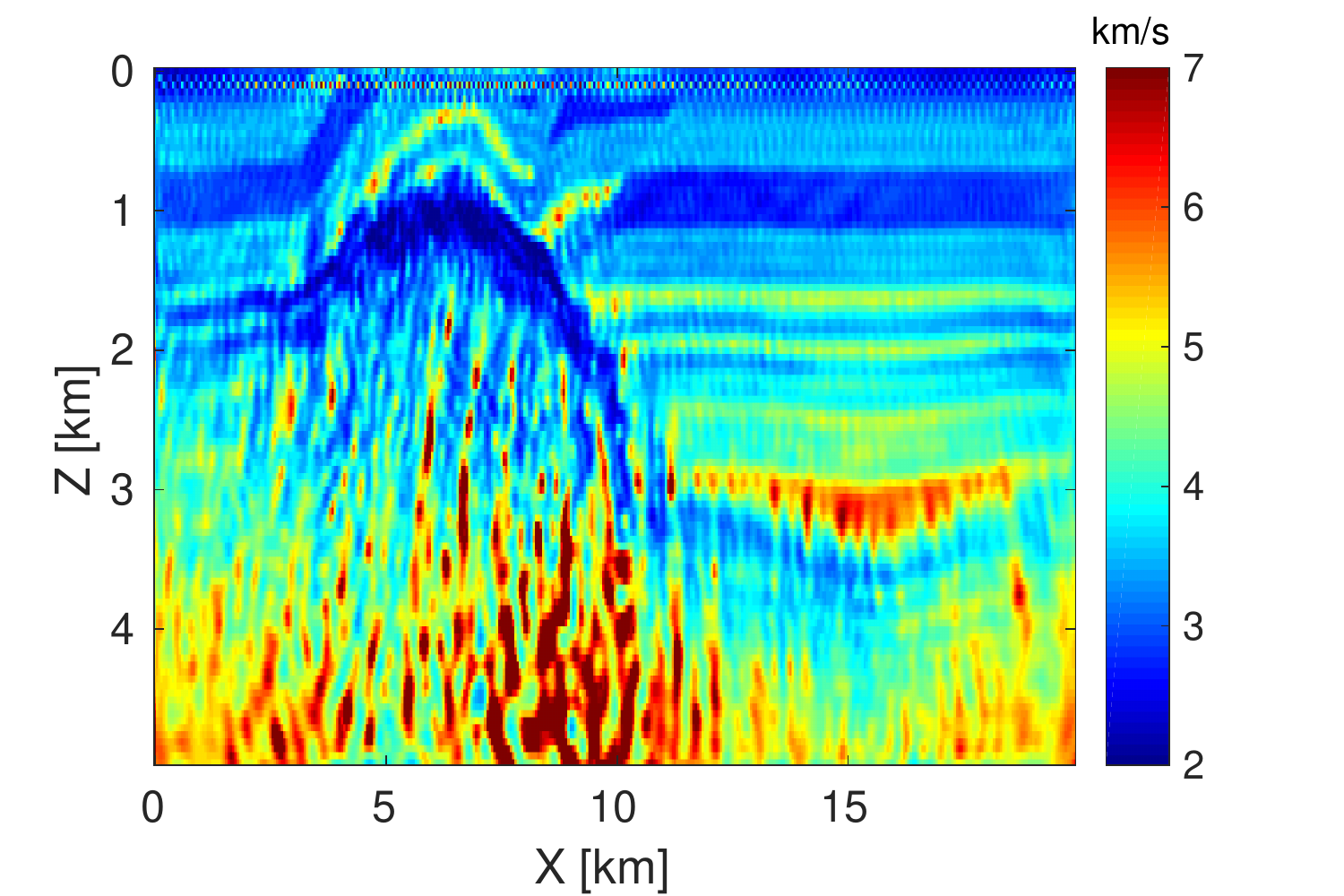}}
\subfloat[Final result of
LRWI\label{fig:ORlrwi03}]{\includegraphics[width=0.300\hsize]{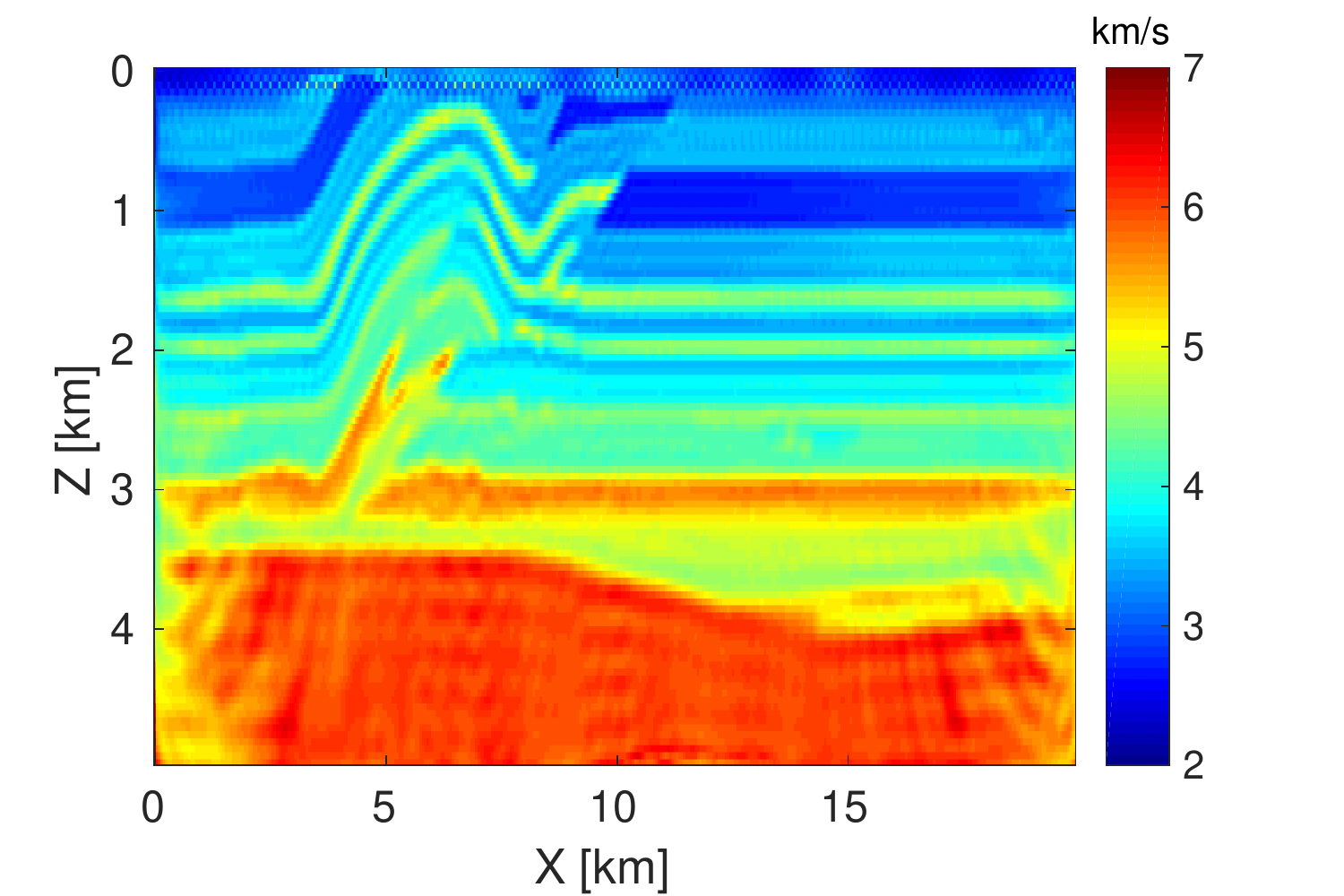}}
\caption{(a) - (c) Results of FWI, WRI, and LRWI after the first
frequency band. (d) - (f) Final results of FWI, WRI, and
LRWI.}\label{fig:OResult0}
\end{figure}

Figure~\ref{fig:OMEpen} shows the comparison of the final relative model
errors for results obtained by conventional WRI using different
$\beta_1$ and LRWI using different $\beta_{1}$ and $\beta_{2}$. We did
not include the result of LRWI with $\beta_2=1e$-16, since the matrix
$\Smatt(\beta_1,\beta_2)^{\top}\Smatt(\beta_1,\beta_2)$ is close to
singular or badly scaled. According to Figure~\ref{fig:OMEpen}, the
selection of $\beta_1 = 1e$-4 produces the best result for WRI, and the
selection of $\beta_1 = 1e$-8 and $\beta_2 = 1e$-12 produces the best
result for LRWI. Figure~\ref{fig:OResult0} shows the results of FWI, WRI
with the best selection of $\beta_1$, and LRWI with the best selection
of $\beta_1$ and $\beta_2$. Figures~\ref{fig:ORfwi01}
to~\ref{fig:ORlrwi01} show the results of the three approaches using the
data of the first frequency band, and Figures~\ref{fig:ORfwi03}
to~\ref{fig:ORlrwi03} show the final results of the three approaches.
Clearly, at the first frequency band, FWI and WRI already converge to
the local minima, while LRWI provides a much better model for the
following inversion, which yields a final result (c.f.
Figure~\ref{fig:ORlrwi03}) that has the minimal relative model error and
matches the true model significantly better than those obtained by FWI
and WRI (c.f. Figures~\ref{fig:ORfwi03} and~\ref{fig:ORwri03}).

Figure~\ref{fig:OMEpen} shows that LRWI obtains inverted models with
relative model errors less than $10\%$ with the selection of
$(\beta_1, \beta_2) = (1e\text{-4}, 1e\text{-12})$ and
$(1e\text{-8}, 1e\text{-12})$, while the relative model errors for
results of WRI are larger than $20\%$. This comparison illustrates that
LRWI with an appropriate relaxation on the rank-1 constraint and the PDE
constraint can mitigate the local minima of FWI and WRI.

\begin{figure}
\centering
\includegraphics[width=0.500\hsize]{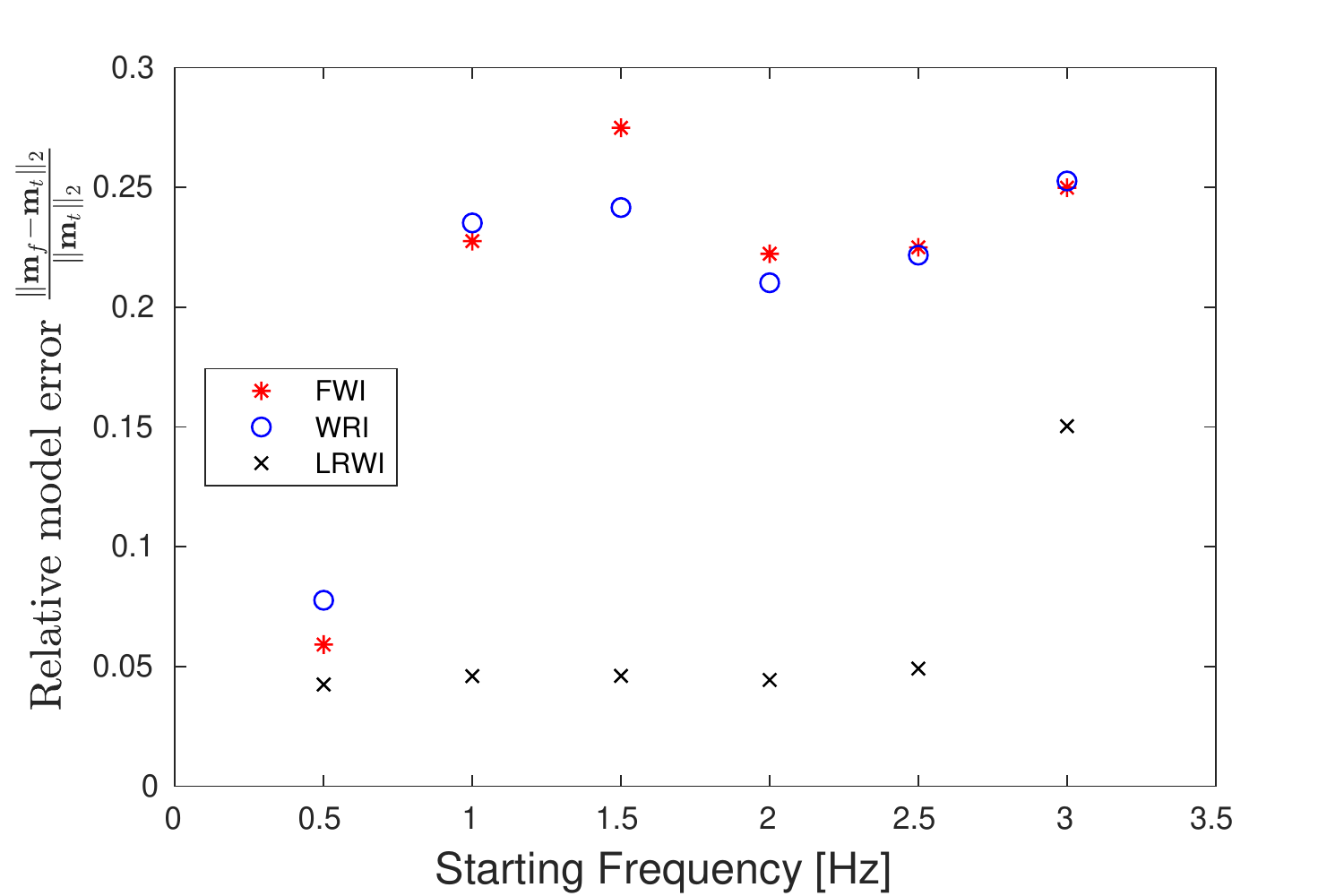}
\caption{Relative model error comparison for FWI($*$), WRI($\circ$), and
LRWI($\times$) using data with different starting
frequencies.}\label{fig:OMFreq}
\end{figure}

\emph{Robustness with respect to the starting frequency} We also conduct
an example to investigate the robustness of the three methods with
respect to the starting frequency for the Overthrust model. In this
example, we vary the starting frequency from 0.5 Hz to 3.0 Hz. We use
the same initial model shown in Figure~\ref{fig:Oini1}.
Figure~\ref{fig:OMFreq} illustrates the relative model errors versus the
starting frequency for all the three methods. The highest starting
frequencies for FWI, WRI, and LRWI to obtain an inverted result with an
acceptable relative model error ($\leq 10\%$) are 0.5 Hz, 0.5 Hz, and
2.5 Hz, respectively. This comparison implies that under the
aforementioned experimental settings, LRWI can conduct a successful
inversion with a starting frequency fifth large as that of conventional
FWI and WRI.

\section{Discussions}\label{discussions}

This paper introduces the basics of a ``Lift'' and ``Relax'' approach
for the waveform inversion with PDE constraints. We have presented
promising initial results in mitigating problems of local minima, while
some aspects of the proposed approach warrant further investigations.

The selection of the penalty parameters $\lambda$ and $\gamma$ are
essential to the success of the proposed LRWI as shown in both numerical
examples. While our analysis and results imply that selecting $\lambda$
to be a small fraction of the largest eigenvalue of
$\Amat^{\top}\Pmat^{\top}\Pmat\Amat^{-1}$ and selecting $\gamma$ to be a
small fraction of the largest fraction between the $\ell_2$-norms of the
vectors $\text{diag}(\Tmat_{i,j}(\lambda))$ and
$\mvect_{i}\odot\mvect_{j}$ at initial iterations yields plausible
results, a more solid justification of this observation would be
desirable for more robust approaches.

While we can use LRWI to conduct inversion with low-frequency data and
produce a good initial model for conventional FWI with a small
computational cost, the application of the LRWI to high-frequency data
can further help us bypass more potential local minima in the objective
function of conventional FWI. To address high-frequency data, a fast
solver for Equation~\ref{AnaSolution} would be worthwhile. One possible
solution is to use efficient direct or iterative solvers designed for
the Helmholtz equation as a preconditioner for the linear system in
Equation~\ref{AnaSolution}.

Finally, the rank of the lifted matrices could be potentially worth
exploring. In this work, we lift the unknown variables from vectors to
rank-2 matrices due to the consideration of storage and computational
cost. Indeed, if the storage and computational cost are not bottlenecks,
we can lift the unknown variables to matrices with higher ranks and
study the effect of the rank on the final inversion results.

\section{Conclusions}\label{conclusions}

We have presented a ``Lift'' and ``Relax'' approach for waveform
inversion problems with PDE constraints. The proposed method is based on
a PDE relaxation and a rank-2 variable relaxation. The reformulation
results in an unconstrained optimization problem with respect to a
rank-2 matrix that contains both lifted model parameters and wavefields.
To avoid storing and updating the rank-2 wavefields during the
optimization, we use the variable projection method to explicitly
eliminate the rank-2 wavefields by solving an overdetermined linear
system. We show that the proposed approach is able to explore a much
larger search space with an acceptable additional computational cost
compared to conventional FWI and WRI.

The main algorithmic difference with conventional FWI and WRI is the
rank-2 variable lifting and the resulting overdetermined system required
to solve. Instead of solving PDEs, we formulate an overdetermined system
of equations that consists of the discretized rank-2 PDE, the
measurements, and the rank-1 regularizations. We study the properties of
this overdetermined system with respect to the selection of the penalty
parameters $\lambda$ and $\gamma$. We show that the condition number of
the overdetermined system can reach a similar value as that of the
original PDE by tuning the two parameters. Therefore, it is plausible
that we can solve the overdetermined system as efficiently using a
similar approach as is applied to the original PDE.

The numerical examples show that the proposed LRWI is able to conduct
successful inversion with higher-frequency data and poorer initial
models compared with conventional FWI and WRI. The numerical examples
further show that through tuning the penalty parameters $\lambda$ and
$\gamma$, the proposed approach can find a search path in the enlarged
space that bypasses the potential local minima in the objective function
of conventional FWI and WRI.

\section{Acknowledgments}\label{acknowledgments}

The authors acknowledge the funding and support provided by ExxonMobil
Research and Engineering Company. Dr.~Laurent Demanet is also supported
by AFOSR grant FA9550-17-1-0316.

\bibliography{zfang}

\end{document}